\documentclass[12pt,reqno,showkeys]{amsart}  
\usepackage{amssymb} 
\usepackage{amsthm}
\usepackage{amsmath}
\usepackage{color}
\allowdisplaybreaks[1]
\newtheorem{thm}{Theorem}[section]
\newtheorem*{def*}{Definition}

\newtheorem{prop}[thm]{Proposition}
\newtheorem{lem}[thm]{Lemma}
\newtheorem{rem}{Remark}[section]

\numberwithin{equation}{section}

\renewcommand{\theequation}{\thesection.\arabic{equation}}

\newenvironment{pr}[1]
   {{\noindent \bf Proof of {#1}.\  }}{\hfill \qed}
                
\newcommand{\pt}{\partial}     
             
\renewcommand{\th}{\theta}                
\newcommand{\D}{\mathcal{D}}
\newcommand{\Ep}{\mathcal E}           
\newcommand{\F}{\mathcal {F}}

\newcommand{\R}{\mathbb R}

\newcommand{\Nt}{\mathbb N}
\newcommand{\N}{\nabla}

\newcommand{\Om}{\Omega}
\newcommand{\al}{\alpha}
\newcommand{\kp}{\kappa}
\newcommand{\gm}{\gamma}
\newcommand{\Gm}{\Gamma}
\renewcommand{\L}{\mathcal{L}}

\newcommand{\ep}{\varepsilon}
\newcommand{\lm}{\lambda}

\newcommand{\del}{\delta}
\newcommand{\Del}{\Delta}

\newcommand{\s}{\sigma}

\renewcommand{\t}{\tau}

\newcommand{\dsp}{\displaystyle}

\newcommand{\bt}{\beta}
\newcommand{\cd}{\cdot}

\newcommand{\E}{\mathcal{E}}

\newcommand{\vp}{\varphi}
\def\<{\langle }

\renewcommand{\qed}{\qquad\kern1pt   
   \vbox{\hrule height 0.6pt      
         \hbox{\vrule width 0.6pt 
               \vbox{\vskip 6pt}  
               \hskip 3pt
              \vrule width 1.3pt} 
         \hrule depth 1.3pt}     
   \kern1pt}

\newcommand{\eqntag}{\addtocounter{equation}{1}\tag{\theequation}} 

\begin{document}
\title{Global existence and boundedness of solutions\\
to a fully parabolic chemotaxis system\\
with indirect signal production in $\mathbb{R}^4$
}
\maketitle

 \vskip5mm
{\small
	\hskip7mm\noindent
	{\normalsize \sf Tatsuya Hosono 
		\footnote{
			ORCID: 0000-0002-1027-9485, Corresponding Author}}
	\vskip1mm\hskip1.5cm
	Mathematical Institute,
	Tohoku University  \hfill\break\indent\hskip1.5cm
	Sendai  980-8578, Japan \hfill\break\indent\hskip1.5cm
	{tatsuya.hosono.q2@dc.tohoku.ac.jp}
}
\vskip5mm
{\small
	\hskip7mm\noindent
	{\normalsize \sf Philippe Lauren\c{c}ot \footnote{ORCID: 0000-0003-3091-8085}}
	\vskip1mm\hskip1.5cm
	Laboratoire de Math\'{e}matiques (LAMA) UMR 5127 \hfill\break\indent\hskip1.5cm 
	Universit\'{e} Savoie Mont Blanc, CNRS \hfill\break\indent\hskip1.5cm
	F--73000 Chamb\'{e}ry, France \hfill\break\indent\hskip1.5cm
	{philippe.laurencot@univ-smb.fr}
	}

\vspace{5mm}
\begin{abstract}
Global existence and boundedness of solutions to the Cauchy problem for
the four dimensional fully parabolic chemotaxis system with indirect signal production are studied. We prove that solutions with initial mass below $(8\pi)^2$ exist globally in time.
This value $(8\pi)^2$ is known as the four dimensional threshold value of the initial mass determining whether blow-up of solutions occurs or not. Furthermore, some condition on the initial mass guaranteeing that the solution remains uniformly bounded is also obtained.
\end{abstract}

\vspace{5mm}

\vspace{5mm}


\vspace{5mm}
\section{Introduction}
\vspace{5mm}
We study the Cauchy problem for the chemotaxis system with indirect signal production:
\begin{equation}
\left\{
\begin{aligned}
&\pt_t u -\Del u +  \N\cd\big(u\N v\big)=0,
& t>0,\, x\in\R^4,
\\
&\pt_t v -d_1\Del v+\lm_1 v =w,
& t>0,\, x\in\R^4,
\\
&\pt_t w -d_2\Del w+\lm_2 w =u,
& t>0,\, x\in\R^4,
\\
&(u,v,w)(0,x)=(u_{0},v_0,w_0)(x),
&\, x\in\R^4,
\end{aligned}
\right.
\label{eqn;FS-model}
\end{equation}
with positive diffusion coefficients $d_1,d_2>0,$ and nonnegative death rates $\lm_1,\lm_2\ge 0$, the initial data $(u_0,v_0,w_0)$ being nonnegative functions on $\R^4$. The purpose of this paper is to show the global existence and the boundedness of nonnegative solutions to the Cauchy problem~\eqref{eqn;FS-model}
under some assumptions on the initial mass of $u_0$.

\vspace{5mm}
Chemotaxis is a biological phenomenon involving migrations, particularly aggregation, 
of cells or organisms in reaction to the concentration gradient of some chemical stimuli.
In 1970, a basic mathematical system modeling chemotaxis was set up by Keller--Segel \cite{Ke-Se} as follows:
\begin{equation}
\left\{
\begin{aligned}
&\pt_t u-\Del u+ \N\cd (u\N v)=0, 
\\
&\pt_t v -\Del v+\lm v=u,
\end{aligned}
\right.
\label{eqn;KS}
\end{equation}
with  $\lm\ge0$. Here, the unknown functions $u$ and $v$ represent the density of the slime molds population and the concentration of the chemical substance, respectively, so that it is natural to assume that the solutions are nonnegative. If we consider the system~\eqref{eqn;KS} in $\R^n$ $(n\ge2)$, the system~\eqref{eqn;KS} with $\lm=0$ is invariant under the scaling transformation
\begin{align*}
u_{\mu}(t,x)=\mu^2 u(\mu^2 t,\mu x),
\quad
v_{\mu}(t,x)=v(\mu^2 t,\mu x),
\quad
\mu>0.
\eqntag
\label{eqn;scaling}
\end{align*} 
According to the Fujita--Kato principle \cite{Ca-We, Fu-Ka}, the scaling critical space of the initial condition $u_0$ of the cell density $u$ is $L^{n/2}(\R^n)$ and it is its size in $L^{n/2}(\R^n)$ which determines the dynamics of the solution. Observing that suitably defined nonnegative solutions to~\eqref{eqn;KS} satisfy the mass conservation $\|u(t)\|_1=\|u_0\|_1$ for $t>0$, we realize that this conservation law matches the critical space $L^{n/2}(\R^n)$ in the two dimensional case~$n=2$ and the behavior of the solutions to~\eqref{eqn;KS} is governed by the size of the initial mass $\|u_0\|_1$ of the first component $u$. Indeed, the problem~\eqref{eqn;KS} in space dimension~$n=2$ features critical phenomena and has been the subject of several studies. In particular, it is known that the solution to~\eqref{eqn;KS} exists for all times if the initial mass $\|u_0\|_1$ is smaller than an explicit threshold value $M^*$, whereas finite (or infinite) time blowup may occur for solutions having an inital mass exceeding $M^*$, see \cite{Bi-Na94, Ca-Co08, Co-Pe-Za, Ga-Za, He-Ve, Mi13, Mi20, Na95, Na-Og11, Na-Se-Yo}.

\vspace{5mm}
More recently, an attraction-repulsion chemotaxis system was proposed in \cite{Lu-Ro-Ke-Mo03} as a mathematical model describing aggregation of microglia cells in the central nervous system during the Alzheimer disease:
\begin{equation}
\left\{
\begin{aligned}
&\pt_t u-\Del u+\N\cd \big(u\N( \bt_1v-\bt_2 w)\big) =0, 
\\
&\pt_t v-\Del v+\lm_1 v=u,
\\
&\pt_t w-\Del w+ \lm_2w=u,
\end{aligned}
\right.
\label{eqn;AR}
\end{equation}
with $\bt_1,\bt_2,\lm_1,\lm_2>0$. The functions $u$, $v$ and $w$ stand for
the density of microglia cells, the concentration of chemoattractant, and the concentration of chemorepellent, respectively. The parameters $\bt_1$ and $\bt_2$ account for the sensitivity of microglia cells to the chemical substances and are known to rule the dynamics of~\eqref{eqn;AR}. Specifically, if $\bt_1<\bt_2$ then global solutions exist whatever the size of the initial data \cite{Ji-Ta, Ji2015, Ya3},
whereas critical phenomena may appear in the two dimensional setting when $\bt_1>\bt_2$, as observed in \cite{Fu-Su,Na-Se-Ya1, Na-Se-Ya2} for instance. It is also worth pointing out here that, in the special case $\bt:=\bt_1=\bt_2$ and $\lm_1\neq \lm_2$, the system~\eqref{eqn;AR} rewrites as follows:
\begin{equation}
\left\{
\begin{aligned}
&\pt_t u -\Del u + \N\cd\big(u\N z\big)=0,
\\
&\pt_t z -\Del z+\lm_1 z= \beta(\lm_2-\lm_1)w,
\\
&\pt_t w -\Del w+ \lm_2w =u,
\end{aligned}
\right.
\label{eqn;FS-model2}
\end{equation}
where $z:=\bt(v-w)$. In particular, the system~\eqref{eqn;FS-model2} exactly corresponds to the system~\eqref{eqn;FS-model} with $d_1=d_2=1$ when $\beta(\lm_2-\lm_1)=1$. A noteworthy difference between~\eqref{eqn;FS-model} (or~\eqref{eqn;FS-model2}) and~\eqref{eqn;KS} is that the density of cells $u$ is indirectly influenced by $v$ (or $z$) through the intermediate species $w$ in the former while the density of cells is directly influenced by the concentration of the chemical substance in the latter. In fact, the system~\eqref{eqn;FS-model} has a quite different structure from~\eqref{eqn;KS} and the invariant scaling space of $u_0$ is given by $L^{n/4}(\R^n)$, since the scaling transformation leaving~\eqref{eqn;FS-model} invariant differs from~\eqref{eqn;scaling} and is now given by 
\begin{align*}
u_{\mu}(t,x)=\mu^4 u(\mu^2 t,\mu x),
\quad
z_{\mu}(t,x)=z(\mu^2 t,\mu x),
\quad
w_{\mu}(t,x)=\mu^2w(\mu^2 t,\mu x),
\quad
\mu>0.
\end{align*}
It is then in the four dimensional setting $n=4$ that the size of  the initial mass $u_0$ of the cell density classifies the global behavior of solutions to~\eqref{eqn;FS-model}.
Fujie--Senba \cite{Fu-Se17, Fu-Se19} consider the initial boundary value problem associated with~\eqref{eqn;FS-model}, supplemented with appropriate boundary conditions, and  investigate the four dimensional critical phenomena for the initial mass. Yamada \cite{Ya1} and Hosono--Ogawa \cite{Ho-Og} (see also \cite{Ho}) study the Cauchy problem for the system~\eqref{eqn;FS-model2} when the second and third equations $z$ and $w$ of~\eqref{eqn;FS-model2} are replaced by their elliptic counterparts and also report such a critical phenomenon. Let us finally point out that lower space dimensions $n\in\{1,2,3\}$ correspond to scaling subcritical cases, so that solutions to~\eqref{eqn;FS-model} exist globally in time and are uniformly bounded, see \cite{Fu-Se17, Ji-Li, Ya2}.

\vspace{5mm}
Besides, a related chemotaxis system is derived in \cite{St-Ty-Po} featuring a non-diffusive species and describing pattern formation of the Mountain Pine Beetles in a forest:
\begin{equation}
\left\{
\begin{aligned}
&\pt_t u -\Del u + \N\cd\big(u\N v\big)=0,
\\
&\pt_t v-\Del v +\del v=w,
\\
&\pt_t w+ \lm w =u,
\end{aligned}
\right.
\label{eqn;MPB}
\end{equation}
where $\del,\lm>0$. The model accounts for the interaction
between flying Mountain Pine Beetles with density $u$, nesting Mountain Pine Beetles with density $w$ and Mountain Pine Beetle pheromones (chemoattractant) $v$, and the evolution of $w$ is given by an ordinary differential equation. The system~\eqref{eqn;MPB} can be regarded as a reduced version of~\eqref{eqn;FS-model} with $d_2=0$ and is actually a chemotaxis system with an indirect signal production mechanism relying on an ordinary differential equation. The existence of a critical mass in two space dimensions for~\eqref{eqn;MPB} is investigated in \cite{Lau19, Ta-Wi17} and in the related result \cite{La-St}, while the higher dimensional structure for~\eqref{eqn;MPB} has not been yet elucidated as far as the authors know.

\vspace{5mm}
In this paper, we deal with the Cauchy problem for the fully parabolic chemotaxis system~\eqref{eqn;FS-model} with an indirect signal production mechanism involving a diffusive partial differential equation and we shall establish the existence of global solutions to~\eqref{eqn;FS-model} in the critical four dimensional case under a mass constraint for $u_0$, and study their boundedness as well. As we shall see below, the study of the system~\eqref{eqn;FS-model} on the whole space requires a different approach compared to that used in a bounded domain in \cite{Fu-Se17, Fu-Se19}, in particular to control the behavior of the solutions as $|x|\to\infty$.

\vspace{5mm}
Let us begin with the definition of mild solutions to~\eqref{eqn;FS-model} in the scaling critical space.

\vspace{5mm}
\begin{def*}[mild solution]
Let $2<p<4$ and $4/3<4p/(p+4)<q<2$.
Given $T\in (0,\infty]$ and  $(u_0,v_0,w_0) \in L^1(\R^4)\times W^{2,2}(\R^4)\times L^2(\R^4)$,
we say that a triplet of functions $(u,v,w)$ defined on $[0,T)\times \R^4$ is a mild solution to~\eqref{eqn;FS-model} on $[0,T)$ if
\begin{enumerate}
	\item $u\in C\big([0,T); L^1(\R^4)\big) \cap C\big((0,T); L^{4/3}(\R^4)\big)$;
	\item $v\in C\big([0,T); W^{2,2}(\R^4)\big)\cap C\big((0,T); W^{2,p}(\R^4)\big)$;
	\item $w\in C\big([0,T); L^2(\R^4)\big) \cap C\big((0,T); W^{1,q}(\R^4)\big)$;
	\item $\dsp\sup_{0<t<T} t^{1/2}\|u(t)\|_{4/3}+\dsp\sup_{0<t<T}t^{1-\frac2p}\|v(t)\|_{W^{2,p}}+\dsp\sup_{0<t<T}t^{\frac32-\frac2q}\|w(t)\|_{W^{1,q}}<\infty$;
	\item $(u,v,w)$ satisfy the integral formulations
	\begin{align*}
	&u(t)=e^{t\Del}u_0-\int_0^t \N\cd e^{(t-s)\Del} \big(u(s)\N v(s)\big) ds,
	\\
	&v(t)= e^{t(d_1\Del-\lm_1)} v_0+\int_0^t e^{(t-s)(d_1\Del-\lm_1)} w(s) ds,
	\\
	&w(t)= e^{t(d_2\Del-\lm_2)} w_0+\int_0^t e^{(t-s)(d_2\Del-\lm_2)} u(s) ds
	\end{align*}
	for $0<t<T$, where $e^{t\Del}$ is the heat semigroup given by
	\begin{align*}
	e^{t\Del} f(x) := G_t*f(x)=\int_{ \R^4}G_t(x-y) f(y)dy,
	\quad
	G_t(x)=\frac{1}{(4\pi t)^{2}} \exp\left(-\frac{|x|^2}{4 t}\right)
	\end{align*}
	for $x\in\R^4$ and $f \in L^1(\R^4)$.
\end{enumerate}
\end{def*}

\vspace{5mm}

The starting point of our analysis is the local well-posedness of~\eqref{eqn;FS-model} in the functional setting described in the previous definition.

\vspace{5mm}

\begin{prop}[local well-posedness]\label{prop;LWP}
	Let $2<p<4$ and $4/3<4p/(p+4)<q<2$.
	For  $(u_0,v_0,w_0)\in L^1(\R^4)\times  W^{2,2}(\R^4) \times  L^2(\R^4)$,
	there exist $T\in (0,\infty]$ and a unique mild solution $(u,v,w)$ to~\eqref{eqn;FS-model} on $[0,T)$ in the sense of the above definition. In addition,
	\begin{equation}
		\text{either $T=\infty$ or $T<\infty$ and } \lim_{t\to T} \left\{ \|u(t)\|_{\frac43} + \|v(t)\|_{W^{2,p}} + \|w(t)\|_{W^{1,q}} \right\} = \infty. \label{eqn;bucrit}
	\end{equation}
	Besides, the solution $(u,v,w)$ has higher regularity and satisfies the system~\eqref{eqn;FS-model} in a classical sense on $(0,T)\times\R^4$.

	Moreover, if $u_0,v_0,w_0 \ge 0$ then $u,v,w\ge 0$ in $(0,T)\times\R^4$ and
	\begin{equation*}
		\|u(t)\|_1=\|u_0\|_1, \qquad t\in [0,T).	\eqntag\label{eqn;mass}
	\end{equation*} 
\end{prop}
\vspace{5mm}

\noindent\textbf{Main theorems.}
Our main theorems read as follows:

\vspace{5mm}
\noindent
\textbf{Existence of global solutions.}

\vspace{5mm}
\begin{thm}\label{thm;GE}
Let $d_1,d_2>0$ and $\lm_1,\lm_2\ge 0$. 
Let $u_0\in L_+^1(\R^4)$, $(1+u_0)\log(1+u_0)\in L^1(\R^4)$, 
and $(v_0,w_0)\in W^{2,2}_+(\R^4)\times L_+^1(\R^4)\cap L^2(\R^4)$, where $X_+$ denotes the positive cone of the Banach lattice $X$.
Suppose that
\begin{align*}
\|u_0\|_1<(8\pi)^2d_1d_2.
\eqntag
\label{eqn;global-assumption}
\end{align*}
Then the problem~\eqref{eqn;FS-model} admits a unique global solution $(u,v,w)$ on $(0,\infty)\times\R^4$.
\end{thm}
\vspace{5mm}

We point out once more that, when $d_1=d_2=1$, the value $(8\pi)^2$ is the four dimensional threshold value of the initial mass $\|u_0\|_1$ of the first component of solutions to~\eqref{eqn;FS-model} above which (finite or infinite)-time blow-up may occur, while such a phenomenon is excluded when $\|u_0\|_1$ lies below this value. The optimal value is derived from the best possible constant $32\pi^2$ of the four dimensional Adams inequality
\begin{align*}
	\sup_{f\in W_0^{2,2}(\Om),\,\, \|(-\Del+I)f\|_2\le1}\int_{ \Om} \left(e^{32\pi^2 f^2} -1\right)dx \le C,
\eqntag
\label{eqn;Adams}
\end{align*}
which is valid for any domain $\Om$ of $\R^4$, see \cite{Ad} and \cite[Theorem~1.4]{Ru-Sa}  for instance. The threshold value $(8\pi)^2$ also comes into play in the analysis of the initial-boundary value problem for the system~\eqref{eqn;FS-model} on a bounded domain with no-flux boundary conditions for $u$ and either homogeneous Neumann or homogeneous Dirichlet boundary conditions for $(v,w)$, as shown in Fujie--Senba \cite{Fu-Se17, Fu-Se19}, still assuming $d_1=d_2=1$. More precisely, on the one hand, they prove that solutions with initial mass $\|u_0\|_1<(8\pi)^2$ exists globally in time by using~\eqref{eqn;Adams}, with an additional restriction to the radially symmetric setting when homogeneous Neumann boundary conditions are prescribed for $(v,w)$ \cite{Fu-Se17}. On the other hand, if the initial mass of $u_0$ is bigger than $(8\pi)^2$, then there exists solutions blowing up in finite or infinite time \cite{Fu-Se19}.

Let us also mention that, for the parabolic-elliptic counterpart of~\eqref{eqn;FS-model}, where the second and third equations are replaced by elliptic equations without the time-derivative term, Yamada \cite{Ya1} and Hosono--Ogawa \cite{Ho-Og} proved that there exists a global solution under $\|u_0\|_1<(8\pi)^2$,
meanwhile if $\|u_0\|_1>(8\pi)^2$ then the solution may blow up in finite time.
No result regarding the finite-time blow-up of the solution to the fully parabolic system~\eqref{eqn;FS-model} seems to have been obtained so far.

\vspace{5mm}
\noindent
\textbf{Strategy of the proof of Theorem~\ref{thm;GE}.}
In the study of the two-dimensional Keller--Segel system~\eqref{eqn;KS},
the combination of the Lyapunov functional and the Trudinger--Moser inequality 
 works well and implies the existence of global solutions as soon as the initial mass is smaller than the threshold value \cite{Ca-Co08, Mi13, Na-Og11}, see also \cite{Ga-Za, Na-Se-Yo}. In contrast, in the present situation, the Adams inequality~\eqref{eqn;Adams} and its variants established in \cite{Fu-Se17, Ru-Sa} do not seem to be the right tool to employ, in particular due to the unboundedness of $\R^4$ and the higher dimensional structure. We thus cannot adapt the approach developed for the two-dimensional Keller--Segel system for showing Theorem~\ref{thm;GE}.

Let us introduce the modified Lyapunov functional $\F$ which is naturally associated with the system~\eqref{eqn;FS-model}:
\begin{align*}
\F(u,v):=\,&
\int_{\R^4}(1+u)\log (1+u) dx
+\frac{1}{2}\int_{\R^4}|\pt_t v|^2dx
+\Ep(v;u),
\eqntag
\label{eqn;lyapunov}
\end{align*}
where $\Ep(v;u)$ is the chemical energy defined as
\begin{align*}
\Ep(v;u):=\frac{d_1d_2}{2}\int_{\R^4}|\Del v|^2dx
+\frac{d_1\lm_2+d_2\lm_1}{2}\int_{\R^4}|\N v|^2 dx
+\frac{\lm_1\lm_2}{2}\int_{\R^4}v^2 dx-\int_{\R^4}uvdx.
\end{align*}
The first term in $\F(u,v)$ is  the modified entropy, which replaces the usual entropy term $\int_{\R^4}u\log u dx$ and has the advantage of being nonnegative. It was introduced by Nagai \cite{Na2011} and has in particular proved useful in the study of the global behavior of solutions to chemotaxis systems on the whole space, see for instance \cite{Ho-Og,La-We-Zh,Mi13,Na2011,Na-Og11,Na-Se-Ya2,Ya3}.
The modified Lyapunov functional $\F$ satisfies the following identity:
\begin{align*}
\F(u(t),v(t)) + \int_0^t\D(u(s),v(s))ds =\F(u_0,v_0)+\frac{1}{4}\int_0^t\int_{\R^4}|\N v(s)|^2 dxds
\eqntag
\label{eqn;lyapunopv-identity}
\end{align*}
for $0<t<T$,
where $\D(u,v)$ is the dissipation term (see Lemma~\ref{lem;lyapunov2} below).
The keypoint in the proof of Theorem~\ref{thm;GE} is to derive the following estimates 
 for the solutions to~\eqref{eqn;FS-model} by using the identity~\eqref{eqn;lyapunopv-identity}  stated in Proposition~\ref{prop;apriori-estimate}: 
\begin{align*}
	\int_{\R^4}(1+u(t))\log(1+u(t))dx+\frac{1}{2}\int_{\R^4}|\pt_t v(t)|^2 dx
+\int_0^t \D(u(s),v(s)) ds
\le C(\tau),
\eqntag
\label{eqn;apriori-solutions}
\end{align*}
for $t\in [0,\tau]\cap [0,T)$, where $C(\tau)$ is a positive constant depending on $\tau$ but not on $T$. As already mentioned, due to the unboundedness of $\R^4$ and the higher space dimension, we cannot apply the arguments which were used either for the initial boundary value problem associated with~\eqref{eqn;FS-model} or for the Cauchy problem for the Keller--Segel system~\eqref{eqn;KS} in order to show estimates~\eqref{eqn;apriori-solutions}. Instead, we introduce the chemical energy minimization as in \cite{Ca-Co08}, see Lemma~\ref{lem;Chemical-EM} below, which we combine with a Brezis--Merle type inequality for the solution of a 4th order elliptic equation (see Hosono--Owaga \cite[Theorem~1.4]{Ho-Og}). The control of the potential energy term $\int_{ \R^4}uvdx$ in $\F(u,v)$ which is included in the chemical energy $\Ep$ is required to derive the estimates~\eqref{eqn;apriori-solutions} and relies on the observation that $\Ep$ is the ``energy'' associated with the 4th order elliptic equation $(-d_1\Del+\lm_1)(-d_2\Del+\lm_2) v_u=u$. In this connection, the chemical energy minimization plays an important role in deriving~\eqref{eqn;apriori-solutions} under the optimal constraint~\eqref{eqn;global-assumption} on the initial mass of $u_0$, as it does for the two dimensional Keller--Segel system in the approach developed in \cite{Ca-Co08}.
It actually allows us to estimate the chemical energy $\Ep(v;u)$ of the solution $(u,v,w)$ to~\eqref{eqn;FS-model} by that involving the solution $v_u$ to the elliptic equation $(-d_1\Del+\lm_1)(-d_2\Del+\lm_2) v_u=u$, i.e.,
\begin{equation*}
\Ep(v;u)\ge \Ep(v_u;u).
\end{equation*}
Together with~\eqref{eqn;lyapunopv-identity}, the Brezis--Merle type inequality, 
and the assumption on the optimal size of the initial mass~\eqref{eqn;global-assumption},
the following bound on the chemical energy is derived:
\begin{equation*}
-\Ep(v(t);u(t)) \le C(\tau), \qquad t\in [0,\tau]\cap [0,T).
\end{equation*}
The estimates~\eqref{eqn;apriori-solutions} then readily follow, since 
\begin{align*}
\int_{\R^4}(1+u(t))\log (1+u(t)) dx
&\, +\frac{1}{2}\int_{\R^4}|\pt_t v(t)|^2dx
+\int_0^t \D(u(s),v(s)) ds
\\
&=\,\F(u(t),v(t))-\Ep(v(t);u(t))+\int_0^t \D(u(s),v(s)) ds \le C(\tau)
\end{align*}
for $t\in [0,\tau]\cap [0,T)$, the $H^1$-norm of $v$ being controlled by parabolic regularity thanks to the $L^1$-bound~\eqref{eqn;mass} on $u$, see Lemma~\ref{lem;vw-bound} and the proof of Proposition~\ref{prop;apriori-estimate}. The next step is to show that
\begin{equation*}
\int_0^t \left[ \| \N \pt_t v(s)\|_2^2 + \|\N w(s)\|_2^2 \right] ds\le C(\tau), \qquad t\in [0,\tau]\cap [0,T).
\end{equation*} 
Owing to the above estimates, $L^p$-estimates for $u$ can be derived for any $p\in (1,\infty)$ by the standard energy method and we are left with using the classical Nash--Moser iteration technique to derive an $L^\infty$-estimate on $u$, after establishing 
a bound on the $L^\infty$-norm for $\N v$.

\vspace{5mm}
\noindent
\textbf{Boundedness of solutions.}

\vspace{5mm}
We note that Theorem~\ref{thm;GE} does not state whether the solutions are uniformly bounded with respect to time or not. A condition for the boundedness of the solutions is given as follows.
\vspace{5mm}
\begin{thm}\label{thm;Bdd}
Let $d_1,d_2>0$ and $\lm_1,\lm_2>0$.
Let $u_0\in L_+^1(\R^4)$, $(1+u_0)\log(1+u_0)\in L^1(\R^4)$, 
and $(v_0,w_0)\in L^1_+(\R^4)\cap W^{2,2}(\R^4)\times L_+^1(\R^4)\cap L^2(\R^4)$.
Suppose that
	\begin{align*}\eqntag
	\|u_0\|_1<\frac{1}{\sqrt{3}}(8\pi)^2d_1d_2.
	\label{eqn;bdd-assumption}
	\end{align*}
	Then the unique solution to the system~\eqref{eqn;FS-model} is bounded on $[0,\infty)$ in all $L^p$-spaces; that is,
\begin{align*}
\sup_{t\ge0}\big(\|u(t)\|_p+\|v(t)\|_p+\|w(t)\|_p\big)
<\infty
\end{align*}
for $1\le p \le \infty$.
\end{thm}

\vspace{5mm}
The mass constraint~\eqref{eqn;bdd-assumption} is clearly stronger than~\eqref{eqn;global-assumption} and we are thus not able to prove that all the global solutions we constructed in Theorem~\ref{thm;GE} are bounded. This contrasts markedly with the initial boundary problem for~\eqref{eqn;FS-model} on a bounded domain with appropriate boundary conditions for which the proof of global existence also provides the boundedness of the solution \cite{Fu-Se17, Fu-Se19}. One of the reasons for this difference is that, on a bounded domain, one can use the entropy $\int_{ \Om}u\log udx$ and the corresponding Lyapunov functional without additional control on the decay at infinity of the initial data and derive directly time-independent bounds. We use instead the modified entropy $\int_{ \Om}(1+u)\log{(1+u)}dx$ to deal with the Cauchy problem~\eqref{eqn;FS-model} in $\R^4$ which only provides time-dependent estimates.
Hence, the analysis of the boundedness for the Cauchy problem is the essential difference from the one for the initial boundary  problem on a bounded domain. Such an issue is already encountered in the analysis of the Cauchy problem for related models, see \cite{Ca-Co08,Mi13,Na-Se-Ya1, Na-Sy-Um, Na-Ya}. We however expect that the sole condition~\eqref{eqn;global-assumption} ensures the boundedness of solutions to~\eqref{eqn;FS-model} in the sense of Theorem~\ref{thm;Bdd}.

\vspace{5mm}
\noindent
\textbf{Strategy of the proof of Theorem~\ref{thm;Bdd}.} At the heart of the proof of Theorem~\ref{thm;Bdd} lies an estimate of the $L^1$-norm of $u (-\Delta)^{-1}u$ in terms of $\|\nabla\sqrt{u}\|_2$ and $\|u\|_1$, from which the mass constraint~\eqref{eqn;bdd-assumption} is computed. Here we combine two functional inequalities, the Hardy--Littlewood--Sobolev and the Sobolev inequalities, with their optimal constants $C_{HLS}$ and $C_S$, (cf. \cite{Gi-Tr,Li-Lo}): 
\begin{align*}
\iint_{\R^4\times\R^4}f(x)|x-y|^{-2}f(y)dxdy\le C_{HLS} \| f \|_{\frac43}^2,
\quad
C_{HLS}=\sqrt{\frac{3}{2}}\pi,
\eqntag
\label{eqn;HLS-inequality}
\end{align*}
and
\begin{align*}
\| f \|_{\frac43} \le C_S \|\N f\|_1,
\quad
C_S=\frac{2^{\frac14}}{4\sqrt{\pi}}.
\eqntag
\label{eqn;Sobolev-inequality-4/3}
\end{align*}
Indeed, it follows from H\"{o}lder's inequality that
\begin{align*}
\int_{\R^4}f (-\Del)^{-1}f dx
=\,&
\frac{1}{4\pi^2}\iint_{\R^4\times\R^4} f(x)|x-y|^{-2}f(y)dxdy
\\
\le\,&  \frac{C_{HLS} C_S^2}{4\pi^2} \|\N f\|_1^2\le  \frac{C_{HLS} C_S^2}{4\pi^2}\|f\|_1\int_{ \R^4}\frac{|\N f|^2}{f}dx,
\eqntag
\label{eqn;bdd-method}
\end{align*}
and 
\begin{align*}
\frac{C_{HLS} C_S^2}{4\pi^2} = \frac{\sqrt{3}}{(8\pi)^2}.
\end{align*}

\begin{rem}\label{rem;optcst}
It is unclear whether the above constant is the largest that we can reached with our approach. We actually conjecture that the constant $C_S^2$ could be improved to $C_{GN}^4/4$, where $C_{GN}$ is the optimal constant of the Gagliardo--Nirenberg inequality \cite{We}
\begin{equation*}
	\| f\|_{\frac83} \le C_{GN} \|f\|_2^{\frac12} \|\N f\|_2^{\frac12}, \qquad f\in W^{1,2}(\R^4).
\end{equation*}
Indeed, the estimate on $\|f\|_{\frac43}^2$ performed in the proof of the inequality~\eqref{eqn;bdd-method} can be replaced by
\begin{align*}
\| f\|_{\frac43}^2=\|f^{\frac12}\|_{\frac83}^4\le C_{GN}^4 \| f^{\frac12}\|_2^2\|\N f^{\frac12}\|_2^2
=\frac{C_{GN}^4}{4}\|f\|_1\int_{ \R^4}\frac{|\N f|^2}{f}dx.
\end{align*}
Nevertheless, as far as the authors know, an explicit value of $C_{GN}$ does not seem to be available.
\end{rem}

\vspace{5mm}
The other key argument is the derivation of a differential inequality for a suitably designed energy functional $\L$ (see Proposition~\ref{prop;bound-L} below) 
which differs from the Lyapunov functional $\F$ defined in~\eqref{eqn;lyapunov}.
We use the smallness condition~\eqref{eqn;bdd-assumption} on the initial mass
and~\eqref{eqn;bdd-method} to show that the dissipative term derived from the modified entropy dominates other terms involving $u$, see the proof of Proposition~\ref{prop;bound-L}. Indeed, by use of the hypothesis~\eqref{eqn;bdd-assumption}, we establish  the following differential inequality for the energy functional $\L$:
\begin{align*}
\frac{d}{dt}\L(u,v,w)+ \del \L(u,v,w) +\del \D_1(u,v,w)\le C
\end{align*}
for some $0<\del<1$ and $C>0$, where $\L(u,v,w)$ and $ \D_1(u,v,w)$ are the energy functional and the associated dissipative term defined in~\eqref{def;L_0} and~\eqref{eqn;dissipative}, respectively. We then conclude that, for all $t\in [0,\infty)$, 
\begin{align*}
\L(u(t),v(t),w(t)) +\int_t^{t+1} \left(\|u(s)\|_{\frac32}^{\frac32}+\| \N \pt_t v(s)\|_2^2+\|\N w(s)\|_2^2\right)ds \le C.
\eqntag
\label{eqn;bound-apriori-vw}
\end{align*}
Once~\eqref{eqn;bound-apriori-vw} is shown, the $L^{3/2}$-bound for $u$ follows from the energy method and the uniform Gronwall inequality, and we next use the Nash--Moser iteration technique to complete the proof of the boundedness of the solutions to~\eqref{eqn;FS-model}.

\vspace{5mm}
\section{Local-in-time  solutions}
\vspace{5mm}

We prove in this section the local well-posedness for the problem~\eqref{eqn;FS-model} in the scaling critical space stated in Proposition~\ref{prop;LWP}.
We fix $(u_0,v_0,w_0)\in L^1(\R^4)\times W^{2,2}(\R^4)\times L^2(\R^4)$ and positive constants $\eta_i$ $(i=1,2,3)$ and $T>0$ to be chosen later depending on $(u_0,v_0,w_0)$. Setting $M:=9(\|u_0\|_1+\|v_0\|_{W^{2,2}}+\|w_0\|_2)$, we also fix 
\begin{equation*}
	2<p<4 \;\;\text{ and }\;\; \frac{4}{3} < \frac{4p}{p+4}<q<2,
\end{equation*}
and define the space $X_T$ by
\begin{equation*}\label{def;Xt}
X_{T} := 
\left\{
\begin{aligned}
&u\in L^{\infty}\big( 0,T ; L^1(\R^4)\big)
\cap C \big( (0,T) ; L^{\frac43} (\R^4)\big),\\
&v\in L^{\infty}\big(0,T; W^{2,2}(\R^4)\big)\cap C \big( (0,T) ; W^{2,p} (\R^4)\big),
\\
&w\in L^{\infty }\big(0,T; L^2(\R^4)\big)\cap C \big( (0,T) ; W^{1,q} (\R^4)\big) ; 
\\
&\sup_{0< t<T}\|u(t)\|_1+
\sup_{0< t<T}\|v(t)\|_{W^{2,2}}+
\sup_{0< t<T}\|w(t)\|_{2}\le M,
\\
&\sup_{0<t<T} t^{\frac12} \|u(t)\|_{\frac43}\le \eta_1,
\,\,
\sup_{0<t<T}t^{1-\frac2p} \| v(t)\|_{W^{2,p}}\le \eta_2,
\,\,
\sup_{0<t<T}t^{\frac32-\frac2q}\|w(t)\|_{W^{1,q}} \le \eta_3
\end{aligned}
\right\}.
\end{equation*}
We note that the space $X_T$ is a complete metric space equipped with the distance 
\begin{align*}
d_{X_T}\left((u,v,w), (\bar{u},\bar{v},\bar{w})\right) 
=\,&\sup_{0<t<T} t^{\frac12} \|u(t)-\bar{u}(t)\|_{\frac43}
+\sup_{0<t<T}t^{1-\frac2p} \| v(t)-\bar{v}(t)\|_{W^{2,p}}
\\
&+\sup_{0<t<T}t^{\frac32-\frac2q}\|w(t)-\bar{w}(t)\|_{W^{1,q}}.
\end{align*}
We next define the mapping $\Phi=(\Phi_1,\Phi_2,\Phi_3):X_T\to X_T$ by
	\begin{align*}
&\Phi_1[u,v,w](t)=e^{t\Del}u_0-\int_0^t \N\cdot e^{(t-s)\Del} \big(u(s)\N v(s)\big) ds,
\\
&\Phi_2[u,v,w](t)= e^{t(d_1\Del-\lm_1)} v_0+\int_0^t e^{(t-s)(d_1\Del-\lm_1)} w(s) ds,
\\
&\Phi_3[u,v,w](t)= e^{t(d_2\Del-\lm_2)} w_0+\int_0^t e^{(t-s)(d_2\Del-\lm_2)} u(s) ds
\end{align*}
for $(u,v,w)\in X_T$ and $t\in (0,T)$.

We first show that $\Phi$ is a contraction mapping from $X_T$ to $X_T$ for a suitable choice of the (free) parameters $(\eta_1,\eta_2,\eta_3,T)$, so as to apply the Banach fixed point theorem.

\vspace{5mm}
\begin{lem}\label{lem;contraction-map}
There are $\eta_1,\eta_2,\eta_3>0$ and $T>0$ small enough such that, for any $(u,v,w)\in X_T$, the following holds:
\begin{align*}
\sup_{0< t<T}\|\Phi_1[u,v,w](t)\|_1+
\sup_{0< t<T}\|\Phi_2[u,v,w](t)\|_{W^{2,2}}+
\sup_{0< t<T}\|\Phi_3[u,v,w](t)\|_{2}\le M,
\eqntag \label{eqn;contraction-M}
\end{align*}
and
\begin{equation*}
	\begin{aligned}
	& \sup_{0<t<T} t^{\frac12} \|\Phi_1[u,v,w](t)\|_{\frac43} \le \eta_1,
	\,\,
	\sup_{0<t<T}t^{1-\frac2p} \|\Phi_2[u,v,w](t)\|_{W^{2,p}}\le \eta_2,
	\\
	& \sup_{0<t<T}t^{\frac32-\frac2q}\|\Phi_3[u,v,w](t)\|_{W^{1,q}} \le \eta_3.
	\end{aligned}
\eqntag \label{eqn;contraction-eta}
\end{equation*}
Besides, for $(u,v,w)\in X_T$ and $(\bar{u},\bar{v},\bar{w})\in X_T$,
\begin{align*}
d_{X_T}\left(\Phi[u,v,w],\Phi[\bar{u},\bar{v},\bar{w}]\right)
\le \frac12 d_{X_T}\left((u,v,w),(\bar{u},\bar{v},\bar{w})\right).
\eqntag\label{eqn;contraction-d}
\end{align*}
\end{lem}
\vspace{5mm}

\begin{pr}{Lemma~\ref{lem;contraction-map}}
Let $(u,v,w)\in X_T$. Owing to H\"{o}lder's inequality, the continuous embedding of $W^{2,2}$ in $W^{1,4}$, and the properties of the heat kernel,
\begin{align*}
\|\Phi_1[u,v,w](t)\|_1
\le\,& \| e^{t\Del} u_0\|_1
+C\int_0^t (t-s)^{-\frac12} \|u(s)\N v(s)\|_1 ds
\\
\le\,& \|  u_0\|_1
+C\int_0^t (t-s)^{-\frac12}\|u(s)\|_{\frac43}\|\N v(s)\|_4 ds
\\
\le\,&\|u_0\|_1+C \int_0^t (t-s)^{-\frac12} s^{-\frac12}ds
\left(\sup_{0< s<t}s^{\frac12}\|u(s)\|_{\frac43}\right)
\left(\sup_{0< s<t}\|v(s)\|_{W^{2,2}}\right)
\\
\le\,&\frac{M}{9}+C_1 \mathsf{B}\left(\frac12,\frac12\right) M \eta_1,
\eqntag \label{eqn;L1}
\end{align*}
where the beta function $\mathsf{B}(x,y)$ for $x,y>0$ is given by
\begin{align*}
\mathsf{B}(x,y):=\int_0^1(1-\t)^{x-1} \t^{y-1}d\t,
\quad
(x,y)\in(0,\infty)^2.
\end{align*}
Choosing $\eta_1>0$ small such that
\begin{equation*}
	C_1\mathsf{B}\left(\frac12,\frac12\right)  \eta_1\le \frac29,
\end{equation*}
we see that
\begin{align*}
\sup_{0< t<T}\|\Phi_1[u,v,w](t)\|_1 \le \frac{M}{3}.
\end{align*}
With regard to the estimate for $\Phi_2[u,v,w]$, the smoothing properties of the heat kernel imply that
\begin{align*}
\|\Phi_2[u,v,w](t)\|_{W^{2,2}}
\le\,&\|v_0\|_{W^{2,2}}+C\int_0^t (t-s)^{-2\left(\frac1q-\frac12\right)-\frac12} s^{-\frac32+\frac2q}ds
\left(\sup_{0< s<t}s^{\frac32-\frac2q}\|w(s)\|_{W^{1,q}}\right)
\\
\le\,& \frac{M}{9} + C_2 \mathsf{B}\left(\frac32-\frac2q, \frac2q-\frac12\right) \eta_3.
\eqntag\label{eqn;W-22}
\end{align*}
Taking $\eta_3>0$ such that
\begin{align*}
C_2 \mathsf{B}\left(\frac32-\frac2q, \frac2q-\frac12\right) \eta_3\le \frac{2}{9}M,
\end{align*}
it follows that
\begin{align*}
\sup_{0< t<T}\|\Phi_2[u,v,w](t)\|_{W^{2,2}} \le \frac{M}{3}.
\end{align*}
Similarly,
\begin{align*}
\|\Phi_3[u,v,w](t)\|_2 \le\,& \|w_0\|_2
+ C \int_0^t (t-s)^{-\frac12} s^{-\frac12} ds \left(\sup_{0< s<t}s^{\frac12}\|u(s)\|_{\frac43}\right)
\\
\le\,&  \frac{M}{9} + C_3 \mathsf{B}\left(\frac12,\frac12\right)\eta_1,
\eqntag\label{eqn;L2}
\end{align*}
so that, assuming further that
\begin{equation*}
	C_3 \mathsf{B}\left(\frac12,\frac12\right)  \eta_1\le \frac29 M,
\end{equation*}
we obtain
\begin{align*}
\sup_{0< t<T}\|\Phi_3[u,v,w](t)\|_2 \le \frac{M}{3}.
\end{align*}
Hence, combining the above estimates provides~\eqref{eqn;contraction-M}.

Next let us show the assertion~\eqref{eqn;contraction-eta}. To this end, let $(r_1,r_2)\in(1,\infty)^2$ be given by
\begin{align*}
\frac{1}{r_1} := \frac34+\frac{1}{r_2},
\quad
\frac1{r_2} := \frac1p-\frac14.
\end{align*}
Then we see that, by H\"{o}lder's inequality and the continuous embedding of $W^{1,p}$ in $L^{r_2}$,
\begin{align*}
t^{\frac12}\left\|\int_0^t \N\cdot e^{(t-s)\Del}(u(s)\N v(s))ds\right\|_{\frac43}
\le\,&Ct^{\frac12}
\int_0^t (t-s)^{-2\left(\frac1{r_1}-\frac34\right)-\frac12} \| u(s)\N v(s)\|_{r_1} ds
\\
\le\,&Ct^{\frac12}
\int_0^t (t-s)^{-\frac2{p}} \|u(s)\|_{\frac43} \|\N v(s)\|_{r_2} ds
\\
\le\,&Ct^{\frac12}
\int_0^t (t-s)^{-\frac2{p}} s^{-\frac12} s^{-1+\frac2p}ds
\\
&\qquad \times \left(\sup_{0< s<t}s^{\frac12}\|u(s)\|_{\frac43}\right)\left(\sup_{0< s<t}s^{1-\frac2p}\|v(s)\|_{W^{2,p}}\right)
\\
\le\,&C_4 \mathsf{B}\left(1-\frac2p,\frac2p-\frac12\right) \eta_2 \left(\sup_{0< s<t}s^{\frac12}\|u(s)\|_{\frac43}\right).
\eqntag \label{eqn;L-4/3}
\end{align*}
Choosing
\begin{equation*}
	C_4 \mathsf{B}\left(1-\frac2p,\frac2p-\frac12\right) \eta_2\le \frac12
\end{equation*}
in~\eqref{eqn;L-4/3} gives
\begin{align*}
t^{\frac12}\|\Phi_1[u,v,w](t)\|_{\frac43}\le
t^{\frac12}\|e^{t\Del}u_0\|_{\frac43} + \frac12 \left(\sup_{0< s<t}s^{\frac12}\|u(s)\|_{\frac43}\right),
\end{align*}
from which we deduce that 
\begin{align*}
\sup_{0< t<\tau}t^{\frac12}\|\Phi_1[(u,v,w)](t)\|_{\frac43}\le
\sup_{0< t<\tau}t^{\frac12}\|e^{t\Del}u_0\|_{\frac43}+\frac12 \left(\sup_{0< t<\tau}t^{\frac12}\|u(t)\|_{\frac43}\right)
\eqntag \label{eqn;L4/3}
\end{align*}
for all $\tau\in (0,T]$. In particular, since
\begin{align*}
\lim_{t\to 0}t^{\frac12}\|e^{t\Del}u_0\|_{\frac43}=0
\end{align*}
according to Weissler \cite[equation~(3.4)]{Weis},
we can choose $T$ small enough such that
\begin{align*}
t^{\frac12}\|e^{t\Del}u_0\|_{\frac43}\le \frac{\eta_1}{2},\quad
t\in(0,T).
\eqntag
\label{eqn;initial-0}
\end{align*}
Combining~\eqref{eqn;L4/3} and~\eqref{eqn;initial-0} leads us to
\begin{equation*}
\sup_{0< t<T}t^{\frac12}\|\Phi_1[u,v,w](t)\|_{\frac43}\le \frac{\eta_1}{2}+\frac{\eta_1}{2}=\eta_1.
\end{equation*}
As for the estimates of $\Phi_2[u,v,w]$ and $\Phi_3[u,v,w]$ for~\eqref{eqn;contraction-eta}, the proof is similar to that of~\eqref{eqn;W-22} and~\eqref{eqn;L2}. Indeed, due to the linearity of $\Phi_2$ and $\Phi_3$,
\begin{align*}
t^{1-\frac2p}\|\Phi_2[u,v,w](t)\|_{W^{2,p}}\le 
t^{1-\frac2p}\|e^{t\Del} v_0\|_{W^{2,p}}
+ C_5 \mathsf{B}\left(\frac12-\frac2q+\frac2p,\frac2q-\frac12\right) \eta_3
\eqntag \label{eqn;W-2p}
\end{align*}
and
\begin{align*}
t^{\frac32-\frac2q}\|\Phi_3[u,v,w](t)\|_{W^{1,q}}
\le\,&t^{\frac32-\frac2q}\|e^{t\Del}w_0\|_{W^{1,q}}
+ C_6 \mathsf{B} \left(\frac2q-1,\frac12\right) \eta_1,
\eqntag \label{eqn;W-1q}
\end{align*}
observing that the right-hand side of~\eqref{eqn;W-2p} is finite due to the constraint $q>4p/(p+4)$. Since
\begin{equation*}
	\lim_{t\to 0}t^{1-\frac2p}\|e^{t\Del} v_0\|_{W^{2,p}}=\lim_{t\to 0}t^{\frac32-\frac2q}\|e^{t\Del}w_0\|_{W^{1,q}}=0,
\end{equation*}
we may choose $\eta_1,\eta_3,T>0$ small enough so that 
\begin{align*}
	C_5 \mathsf{B}\left(\frac12-\frac2q+\frac2p,\frac2q-\frac12\right) \eta_3 \le \frac{\eta_2}{2}, \quad C_6 \mathsf{B} \left(\frac2q-1,\frac12\right) \eta_1 \le \frac{\eta_3}{2},
\end{align*}
and 
\begin{equation*}
	\sup_{0< t<T}t^{1-\frac2p}\|e^{t\Del} v_0\|_{W^{2,p}} \le \frac{\eta_2}{2}, \quad \sup_{0<t<T}t^{\frac32-\frac2q}\|e^{t\Del}w_0\|_{W^{1,q}} \le \frac{\eta_3}{2},
\end{equation*}
and thereby deduce from~\eqref{eqn;W-2p} and~\eqref{eqn;W-1q} that 
\begin{equation*}
	\sup_{0<t<T}t^{1-\frac2p} \|\Phi_2[u,v,w](t)\|_{W^{2,p}}\le \eta_2,
	\quad
	\sup_{0<t<T}t^{\frac32-\frac2q}\|\Phi_3[u,v,w](t)\|_{W^{1,q}} \le \eta_3,
\end{equation*}
which  complete the proof of the assertion~\eqref{eqn;contraction-eta}.

Finally, an analogous argument implies the assertion~\eqref{eqn;contraction-d}.
\end{pr}
\vspace{5mm}

\begin{lem}\label{lem;local-wellposed}
There is a unique mild solution $(u,v,w)$ to the Cauchy problem~\eqref{eqn;FS-model} defined on a maximal time interval $[0,T)$ with $T\in (0,\infty]$ which satisfies the alternative~\eqref{eqn;bucrit} and depends continuously on the inital data. More precisely, given $(u_0,v_0,w_0)$ and $(\bar{u}_0, \bar{v}_0, \bar{w}_0)$ in $L^1(\R^4)\times W^{2,2}(\R^4)\times L^2(\R^4)$, we denote by $(u,v,w)$ and $(\bar{u},\bar{v},\bar{w})$ the corresponding solutions to~\eqref{eqn;AR}, respectively, with respective maximal existence time $T$ and $\bar{T}$. For any $0<t_0<\min\{T,\bar{T}\}$, there is $C(t_0)>0$ such that
\begin{align*}
&\sup_{0< t<t_0}\| u(t)-\bar{u}(t)\|_1 \le C(t_0) \left(\|u_0-\bar{u}_0\|_1+\|v_0-\bar{v}_0\|_{W^{2,2}}+ \|w_0-\bar{w}_0\|_2\right),
\\
&\sup_{0< t<t_0}\| v(t)-\bar{v}(t)\|_{W^{2,2}} \le C(t_0) \left(\|u_0-\bar{u}_0\|_1+\|v_0-\bar{v}_0\|_{W^{2,2}}+ \|w_0-\bar{w}_0\|_2\right),
\\
&\sup_{0< t<t_0}\| w(t)-\bar{w}(t)\|_2 \le C(t_0) \left(\|u_0-\bar{u}_0\|_1+\|v_0-\bar{v}_0\|_{W^{2,2}}+ \|w_0-\bar{w}_0\|_2\right).
\end{align*}
\end{lem}

\vspace{5mm}
\begin{pr}{Lemma~\ref{lem;local-wellposed}}
By Lemma~\ref{lem;contraction-map}, there are $\eta_1$, $\eta_2$, $\eta_3>0$ and $T>0$ small enough such that $\Phi$ is a contraction mapping from $X_T$ to $X_T$. Hence one can apply the Banach fixed point theorem and conclude that $\Phi$ has a unique fixed point $(u,v,w)\in X_T$ such that
\begin{align*}
\Phi_1[u,v,w]=u,
\quad
\Phi_2[u,v,w]=v,
\quad
\Phi_3[u,v,w]=w.
\end{align*}
In particular, $\Phi_1[u,v,w]=u$ and~\eqref{eqn;L4/3} implies that
\begin{align*}
\sup_{0< s<t} s^{\frac12}\|u(s)\|_{\frac43}\le 2 \sup_{0< s<t} s^{\frac12}\|e^{s\Del}u_0\|_{\frac43}
\end{align*}
for $0<t<T$. Let us show that $u$ belongs to $C\big([0,T); L^1(\R^4)\big)$. For any $\ep>0$, there exists $t_\ep\in (0,T)$ such that
\begin{align*}
\sup_{0< t<t_\ep} t^{\frac12}\| e^{t\Del} u_0\|_{\frac43}\le \frac{\ep}{2},
\end{align*}
from which it follows, by the above estimate, that
\begin{align*}
\sup_{0< s<t_\ep} s^{\frac12}\|u(s)\|_{\frac43}\le 2 \sup_{0< s<t_\ep} s^{\frac12}\|e^{s\Del}u_0\|_{\frac43}
\le \ep.
\end{align*}
Consequently,
\begin{align*}
\lim_{t\to0}\sup_{0< s<t}s^{\frac12}\|u(s)\|_{\frac43} = 0.
\eqntag\label{eqn;L4/3-ep-0}
\end{align*}
Arguing as in the proof of~\eqref{eqn;L1}, we see that
\begin{align*}
\|u(t)-u_0\|_1\le \|e^{t\Del}u_0-u_0\|_1+C \left(\sup_{0< s<t}\|v(s)\|_{W^{2,2}}\right) \left(\sup_{0< s<t}s^{\frac12}\|u(s)\|_{\frac43}\right)
\end{align*}
and it follows from~\eqref{eqn;L4/3-ep-0} and the continuity properties of the heat semigroup that
\begin{align*}
\lim_{t\to0}\|u(t)-u_0\|_1=0.
\end{align*}
Therefore $u\in C\big([0,T); L^1(\R^4)\big)$. Using~\eqref{eqn;L4/3-ep-0},we analogously show from~\eqref{eqn;W-1q} that
\begin{equation}
	\lim_{t\to0}\sup_{0< s<t}s^{\frac32 - \frac{2}{q}} \|w(s)\|_{W^{1,q}} = 0.
	\label{eqn;W1q-ep-0}
\end{equation}
Thanks to~\eqref{eqn;L4/3-ep-0} and~\eqref{eqn;W1q-ep-0}, we may argue as in the proofs of~\eqref{eqn;W-22} and~\eqref{eqn;L2} to establish that 
$v\in C\big([0,T); W^{2,2}(\R^4)\big)$ and $w\in C\big([0,T); L^2(\R^4)\big)$, thereby showing that $(u,v,w)$ is a mild solution to~\eqref{eqn;FS-model} on $[0,T]$. 

We next show the uniqueness of solutions. Let $(u,v,w)$ and $(\bar{u},\bar{v},\bar{w})$ be two solutions to~\eqref{eqn;FS-model} on $[0,T]$ corresponding to the initial data $(u_0,v_0,w_0)$. Noting that the mild formulation provides the representation formulas for $v$ and $w$
\begin{equation*}
v(t)-\bar{v}(t)=\int_0^t e^{(t-s)(d_1\Del-\lm_1)} (w(s)-\bar{w}(s))ds,
\end{equation*}
and
\begin{equation*}
w(t)-\bar{w}(t)=\int_0^t e^{(t-s)(d_2\Del-\lm_2)} (u(s)-\bar{u}(s))ds,
\end{equation*}
a similar argument to~\eqref{eqn;W-2p} and~\eqref{eqn;W-1q} implies that
\begin{align*}
\sup_{0< s<t}s^{1-\frac2p}\| v(s) - \bar{v} (s) \|_{W^{2,p}} 
\le\,& C \sup_{0< s<t} s^{\frac32-\frac2q} \| w(s) -\bar{w}(s) \|_{W^{1,q}}
\\
\le\,& C \sup_{0< s<t} s^{\frac12} \| u(s) -\bar{u}(s) \|_{\frac43},
\end{align*}
so that, using now the representation formula for $u-\bar{u}$ as in the proof of~\eqref{eqn;L-4/3},
\begin{align*}
\sup_{0< s<t}s^{\frac12}\|u(s)-\bar{u}(s)\|_{\frac43}
\le\,&C \left(\sup_{0< s<t}s^{1-\frac2p}\|v(s)\|_{W^{2,p}}\right)
\left(\sup_{0< s<t}s^{\frac12}\|u(s)-\bar{u}(s)\|_{\frac43}\right)
\\
&+ C\left(\sup_{0< s<t}s^{\frac12}\|u(s)\|_{\frac43}\right)\left(\sup_{0< s<t}s^{1-\frac2p}\|v(s)-\bar{v}(s)\|_{W^{2,p}}\right)
\\
\le\,&C \left( \sup_{0< s<t}s^{\frac12}\|u(s)\|_{\frac43}+\sup_{0< s< t}s^{1-\frac2p}\|v(s)\|_{W^{2,p}} \right)
 \left(\sup_{0< s<t}s^{\frac12}\|u(s)-\bar{u}(s)\|_{\frac43}\right).
\end{align*}
Since there exists $t_0\in(0,T)$ sufficiently small such that
\begin{align*}
C\left(\sup_{0< t< t_0}t^{\frac12}\|u(t)\|_{\frac43}+\sup_{0< t< t_0}t^{1-\frac2p}\|v(t)\|_{W^{2,p}} \right)\le \frac12,
\end{align*}
we have
\begin{align*}
\sup_{0< t<t_0} t^{\frac12} \| u(t)-\bar{u}(t)\|_{\frac43}=0,
\end{align*}
from which we conclude that $u(t)=\bar{u}(t)$ for $0<t<t_0$, which in turn implies that $(v,w)(t)=(\bar{v},\bar{w})(t)$ for $0<t<t_0$. Define $t_1\ge 0$ by
\begin{equation*}
	t_1 := \sup\left\{ \tau>0\,;\,(u,v,w)(t)=(\bar{u},\bar{v},\bar{w})(t),\,\text{for all  } 0\le t<\tau \right\}.
\end{equation*}
Clearly $0<t_0\le t_1\le T$. Let us then assume for contradiction that $t_1<T$. Repeating the above argument with initial data $(u,v,w)(t_1)=(\bar{u},\bar{v},\bar{w})(t_1)$, there exists $t_2\in(0,T-t_1)$ such that $(u,v,w)(t)=(\bar{u},\bar{v},\bar{w})(t)$ for all $0 \le t \le t_1+t_2$, which contradicts the definition of $t_1$. Consequently, $(u,v,w)(t)=(\bar{u},\bar{v},\bar{w})(t)$ for $0<t<T$ and we have established the claimed uniqueness.

The continuous dependence of the initial data also follows in the same way. Namely,
for any $(u_0,v_0,w_0)$ and $(\bar{u}_0, \bar{v}_0, \bar{w}_0)$ in $L^1(\R^4)\times W^{2,2}(\R^4)\times L^2(\R^4)$, let $(u,v,w)$ and $(\bar{u},\bar{v},\bar{w})$ be the corresponding solutions to~\eqref{eqn;AR}, respectively, with respective maximal existence time $T$ and $\bar{T}$. For $0<t_0<\min\{T,\bar{T}\}$ and $0<t<t_0$, a similar argument as above follows that
\begin{align*}
&\sup_{0< s<t} s^{\frac12} \| u(s) - \bar{u}(s)\|_{\frac43} 
\\
\le\,& C \| u_0 - \bar{u}_0\|_1+C\left(\sup_{0< s<t}s^{1-\frac2p}\|v(s)\|_{W^{2,p}}\right)
\left(\sup_{0< s<t}s^{\frac12}\|u(s)-\bar{u}(s)\|_{\frac43}\right)
\\
&+C \left(\sup_{0< s<t}s^{\frac12}\|u(s)\|_{\frac43}\right)\left(\sup_{0< s<t}s^{1-\frac2p}\|v(s)-\bar{v}(s)\|_{W^{2,p}}\right)
\\
\le\,&C \| u_0 - \bar{u}_0\|_1+C\left(\sup_{0< s<t}s^{1-\frac2p}\|v(s)\|_{W^{2,p}}\right)
\left(\sup_{0< s<t}s^{\frac12}\|u(s)-\bar{u}(s)\|_{\frac43}\right)
\\
&+C \left(\sup_{0< s<t}s^{\frac12}\|u(s)\|_{\frac43}\right)
\left(\|v_0-\bar{v}_0\|_{W^{2,2}}+\|w_0-\bar{w}_0\|_2+\sup_{0< s<t}s^{\frac12}\|u(s)-\bar{u}(s)\|_{\frac43}\right),
\end{align*}
which implies that
\begin{align*}
\sup_{0< s<t}s^{\frac12}\|u(s)-\bar{u}(s)\|_{\frac43}\le
C(t_0) \left(\|u_0-\bar{u}_0\|_1+\|v_0-\bar{v}_0\|_{W^{2,2}}+ \|w_0-\bar{w}_0\|_2\right).
\end{align*}
Moreover, we observe
\begin{align*}
\sup_{0< s<t}\| v(t) - \bar{v}(s)\|_{W^{2,2}}
\le\,& C \| v_0 - \bar{v}_0\|_{W^{2,2}} + C \sup_{0< s<t} s^{\frac32-\frac2q}\|w(s)-\bar{w}(s)\|_{W^{1,q}}
\\
\le\,& C \| v_0 - \bar{v}_0\|_{W^{2,2}} + C \|w_0 - \bar{w}_0\|_2+C\sup_{0< s<t} s^{\frac12} \| u(s) - \bar{u}(s)\|_{\frac43} 
\\
\le\,&C \left(\|u_0-\bar{u}_0\|_1+\|v_0-\bar{v}_0\|_{W^{2,2}}+ \|w_0-\bar{w}_0\|_2\right).
\end{align*}
Hence,
\begin{align*}
\sup_{0< s<t}\| u(s)-\bar{u}(s)\|_1 
\le\,& C \|u_0-\bar{u}_0\|_1+C\left(\sup_{0< s<t}s^{\frac12}\|u(s)-\bar{u}(s)\|_{\frac43}\right)
\left(\sup_{0< s<t}\|v(s)\|_{W^{2,2}}\right)
\\
&+C\left(\sup_{0< s<t}s^{\frac12}\|u(s)\|_{\frac43}\right)
\left(\sup_{0< s<t}\|v(s)-\bar{v}(s)\|_{W^{2,2}}\right)
\\
\le\,&C(t_0) \left(\|u_0-\bar{u}_0\|_1+\|v_0-\bar{v}_0\|_{W^{2,2}}+ \|w_0-\bar{w}_0\|_2\right),
\end{align*}
and similarly
\begin{align*}
&\sup_{0< s<t}\| w(s)-\bar{w}(s)\|_2 \le C(t_0) \left(\|u_0-\bar{u}_0\|_1+\|v_0-\bar{v}_0\|_{W^{2,2}}+ \|w_0-\bar{w}_0\|_2\right).
\end{align*}
As a consequence, the local well-posedness for~\eqref{eqn;FS-model} follows.

Finally, classical arguments based on the proof of Lemma~\ref{lem;contraction-map} lead to the blowup criterion~\eqref{eqn;bucrit}.
\end{pr}

\vspace{5mm}

\begin{pr}{Proposition~\ref{prop;LWP}}
By virtue of Lemmas~\ref{lem;contraction-map} and~\ref{lem;local-wellposed},
we have established the existence of a unique mild solution to~\eqref{eqn;FS-model} defined on a maximal time interval $[0,T)$ with $T\in (0,\infty]$. As for the regularity of the solutions, we use the standard iteration argument with respect to the derivative. Define $|\N|^\al f(x):=\F^{-1}[ \xi\mapsto |\xi|^\al \F f(\xi)\,](x)$  for $x\in\R^4$ and $\al>0$, where $\F$ denotes the Fourier transform. Let $\tau\in (0,T)$ and $t\in (0,\tau)$. Recalling that, for $2<p<4$,
\begin{equation*}
	U_0(\tau) := \sup_{0< s<\tau}s^{\frac12}\|u(s)\|_{\frac43}<+\infty, \quad V_{p,0}(\tau) := \sup_{0< s<\tau}s^{1-\frac2p}\|v(s)\|_{W^{2,p}}<+\infty,
\end{equation*}
by Lemma~\ref{lem;contraction-map}, we may argue as in the proof of~\eqref{eqn;L-4/3} to obtain
\begin{align*}
&t^{\frac{1+\al}{2}} \| |\N|^{\al} u(t) \|_{\frac43}
\\
\le\,& t^{\frac{1+\al}{2}}  \| |\N|^\al e^{t\Del}u_0\|_{\frac43} + C U_0(\tau) V_{p,0}(\tau) t^{\frac{1+\al}{2}} \int_0^t  (t-s)^{-\frac2{p}-\frac\al2} s^{-\frac{1}2} s^{-1+\frac2p}ds
\\
\le\,&C\|u_0\|_1 + C  U_0(\tau) V_{p,0}(\tau) \mathsf{B}\left(1-\frac2p-\frac\al2,\frac2p-\frac12\right),
\end{align*}
where $2<p<4$ and
\begin{align*}
\al <2-\frac4p<1,
\end{align*}
which implies that
\begin{align*}
U_\al(\tau):=\sup_{0< t<\tau} t^{\frac{1+\al}{2}} \| |\N|^{\al} u(t)\|_{\frac43} <+\infty,
\quad
\al<2-\frac4p.
\end{align*}
This property yields additional regularity estimates for $w$. Indeed, for $\bt\ge \al$, 
\begin{align*}
\| |\N|^{\bt} w(t)\|_{W^{1,q}}
\le\,& C t^{-\frac32+\frac2q-\frac\bt2} \|w_0\|_2
+C U_\al(\tau) \int_0^t (t-s)^{-2\left(\frac34-\frac1q\right)-\frac{1+\bt-\al}{2}} s^{-\frac{1+\al}{2}}ds \\
=\,& C t^{-\frac32+\frac2q-\frac\bt2} \|w_0\|_2
+ C U_\al(\tau) t^{-\frac32+\frac2q-\frac\bt2} \mathsf{B}\left(\frac2q-\frac{\bt-\al}{2}-1, \frac{1-\al}{2}\right),
\end{align*}
where
\begin{align*}
\frac{4}{3} < \frac{4p}{p+4}<q<2, \quad
\al\le \bt<\al +\frac4q-2,
\end{align*}
so that
\begin{align*}
W_{q,\bt}(\tau) := \sup_{0< t<\tau} t^{\frac32-\frac2q+\frac{\bt}{2}} \| |\N|^{\bt} w(t)\|_{W^{1,q}} <+\infty,
\quad \al\le \bt < \al +\frac4q-2.
\end{align*}
Similarly, for $\gm\ge \beta$,
\begin{align*}
\| |\N|^{\gm} v(t) \|_{W^{2,p}}
\le\,&C t^{-1+\frac2p-\frac\gm2} \| v_0\|_{W^{2,2}} + C W_{q,\bt}(\tau) \int_0^t (t-s)^{-2\left(\frac1q-\frac1p\right)-\frac{1+\gm-\bt}{2} } s^{-\frac32+\frac2q-\frac\bt2}ds 
\\
=\,&C t^{-1+\frac2p-\frac\gm2} \| v_0\|_{W^{2,2}}
+C W_{q,\bt}(\tau)  t^{-1+\frac2p-\frac\gm2} \mathsf{B}\left(\frac2p-\frac2q+\frac{1-\gm+\bt}{2}, \frac2q-\frac{\bt+1}{2}\right),
\end{align*}
provided that
\begin{align*}
\beta\le \gm< \bt+\frac4p-\frac4q+1,
\quad
\bt<\al+\frac4q-2<\frac4q-1,
\end{align*}
from which it follows that 
\begin{align*}
V_{p,\gm}(\tau) := \sup_{0< t<\tau} t^{1-\frac2p+\frac\gm2} \| |\N|^{\gm} v(t) \|_{W^{2,p}}<+\infty,
\quad \bt \le \gm<\frac4p-\frac4q+\bt+1 \in (\beta,1).
\end{align*}
The previous analysis leads us to a higher regularity estimate for $u$. Indeed, let $\th\ge\al$. Since $\al\le\bt\le\gm$,
\begin{align*}
& \| |\N|^{\th}u(t)\|_{\frac43} \\
\le\,&C t^{-2\left(1-\frac34\right)-\frac\th2}\|u_0\|_1
\\
&+C\int_0^t(t-s)^{-\frac2p-\frac{\th-\al}{2}} \left(\||\N|^{\al} u(s)\|_{\frac43}\|v(s)\|_{W^{2,p}}+\|u(s)\|_{\frac43}\| |\N|^{\al} v(s)\|_{W^{2,p}}\right)ds
\\
\le\,&C t^{-\frac{1+\th}{2}}\|u_0\|_1
+ C \big(U_0(\tau) V_{p,\al}(\tau) + U_\al(\tau) V_{p,0}(\tau)\big) \int_0^t(t-s)^{-\frac2p-\frac{\th-\al}{2}} s^{-\frac32+\frac2p-\frac{\al}{2}}ds
\\
=\,&C t^{-\frac{1+\th}{2}}\|u_0\|_1 + C \big(U_0(\tau) V_{p,\al}(\tau) + U_\al(\tau) V_{p,0}(\tau)\big) t^{-\frac{1+\th}{2}} \mathsf{B}\left(1-\frac2p-\frac{\th-\s}{2},-\frac12+\frac2p-\frac\s2\right),
\end{align*}
where
\begin{align*}
\al \le \th<2-\frac4p+\al<2,\quad
\al<\frac4p-1.
\end{align*}
Hence, 
\begin{align*}
U_{\th}(\tau) = \sup_{0< t< \tau}t^{\frac{1+\th}{2}}\||\N|^\th u(t)\|_{\frac43}<+\infty,
\quad
\al \le \th<2-\frac4p+\al<2,\quad
\al<\frac4p-1.
\end{align*}
Iterating further the above arguments, we obtain that the solution possesses higher regularity. In fact, using the uniqueness result established in Lemma~\ref{lem;local-wellposed}, we may solve the problem~\eqref{eqn;FS-model} with the initial data $(u(t_0),v(t_0),w(t_0)))$ for some $t_0\in (0,T)$. We thus eventually show that for $s>0$, $u(t)\in L^1(\R^4)\cap W^{s,4/3}(\R^4)$,
$v(t)\in  W^{2,2}(\R^4)\cap W^{s,p}(\R^4)$
and
$w(t)\in L^2(\R^4)\cap W^{s,q}(\R^4)$
 for $t_0<t<T$, 
 and satisfies the system~\eqref{eqn;FS-model} in a  classical sense on $(0,T)\times\R^4$.
 
Let  $(u_0,v_0,w_0)$ be nonnegative initial conditions and
$\{u_{0,k}\}_{k\in\Nt }$,
$\{v_{0,k}\}_{k\in\Nt }$, and $\{w_{0,k}\}_{k\in\Nt }$ be sequences
of nonnegative functions in $C_0^{\infty}(\R^n)$ converging to $u_0$ in $L^{1}(\R^4)$, $v_0$ in $W^{2,2}(\R^4)$, and
$w_0$ in $L^2(\R^4)$ as $k\to\infty$, respectively. Then, for each $k\in\Nt$, we can construct the unique, smooth and integrable solution
to~\eqref{eqn;AR} corresponding to the initial data $(u_{0,k},v_{0,k},w_{0,k})$.
By the parabolic regularity theory,
the solutions $(u_k,v_k,w_k)$ belong to $C^{\infty}( [0,T)\times \R^4)$
and are nonnegative on $(0,T)\times\R^4$. In addition, by the continuous dependence on the initial data established in Lemma~\ref{lem;local-wellposed}, we have the convergences $u_k\to u$ in $C\big([0,T); L^1(\R^4)\big)$, $v_k\to v$ in $C\big([0,T); W^{2,2}(\R^4)\big)$, and $w_k\to w$ in $C\big([0,T); L^{2}(\R^4)\big)$ as $ k\to \infty$,
where $(u,v,w)$ is the solution to~\eqref{eqn;FS-model} with initial data $(u_0,v_0,w_0)$. Hence, one can find subsequences of the approximated solutions
which converge to $(u,v,w)$ almost everywhere in $(0,T)\times\R^4$, which implies that $u$,$v$, and $w$ are nonnegative on $(0,T)\times\R^4$.
Moreover the mass conservation law also holds by integrating $u$ over $\R^4$.
\end{pr}


\vspace{5mm}
We next exploit the fully parabolic structure of~\eqref{eqn;FS-model}, along with the conservation of mass~\eqref{eqn;mass}, to derive a first set of estimates on solutions to~\eqref{eqn;FS-model} which do not depend on their maximal existence time.
\vspace{5mm}
\begin{lem}\label{lem;vw-bound}
Let $(u,v,w)$ be the solution to~\eqref{eqn;FS-model} on $[0,T)$ corresponding to the initial data $(u_0,v_0,w_0)\in L^1(\R^4)\times  W^{2,2}(\R^4) \times L^1(\R^4)\cap L^2(\R^4)$. Then, for $1\le p <2$, there exists a positive constant $C(p)$ such that
\begin{align*}
\|w(t)\|_p \le C(p)(\|w_0\|_p+\|u_0\|_1\mathcal{I}_{\lm_2,p}(t)), \qquad t\in [0,T),
\end{align*}
where 
\begin{equation*}
\mathcal{I}_{\lm_2,p}(t):=\int_0^t e^{-\lm_2 s} s^{-2\left(1-\frac1p\right)}  ds
\le
\left\{
\begin{aligned}
&C(p) \Gm\left(\frac2p-1\right)
&\text{if}\quad
&\lm_2>0,
\\
&C(p) t^{\frac2p-1}
&\text{if}\quad
&\lm_2=0.
\end{aligned}
\right.
\end{equation*}
Moreover, for $2\le q <\infty$,  there exists a positive constant $C(q)$ such that
\begin{equation*}
\|v(t)\|_{q} \le C(q) \left(\|v_0\|_q+\|w_0\|_{\frac{2q}{q+1}}\,\mathcal{I}_{\lm_1,q}(t)
+\|u_0\|_1 \,\mathcal{I}_{\lm_1,q}(t)\,\mathcal{I}_{\lm_2,\frac{2q}{q+1}}(t) \right), \qquad t\in [0,T),
\end{equation*}
where
\begin{equation*}
\mathcal{I}_{\lm_1,q}(t):= \int_0^t e^{-\lm_1s} s^{-1+\frac1q}ds
\le 
\left\{
\begin{aligned}
&C(q) \Gm\left(\frac1q\right)
&\text{if}\quad
&\lm_1>0,
\\
&C(q) t^{\frac1q}
&\text{if}\quad
&\lm_1=0.
\end{aligned}
\right.
\end{equation*}
Besides,
\begin{align*}
\|\N v(t)\|_{2} \le C(\| \N v_0\|_2 + \|w_0\|_{\frac32} \mathcal{I}_{\lm_1,6}(t) + \|u_0\|_1 \mathcal{I}_{\lm_1,6}(t) \mathcal{I}_{\lm_2,\frac{3}{3}}(t)).
\end{align*}
\end{lem}

\vspace{5mm}
\begin{pr}{Lemma~\ref{lem;vw-bound}}
Using the $L^1$-$L^p$ estimate for the heat semigroup and~\eqref{eqn;mass}, the following estimate holds:
\begin{align*}
\left\|\int_0^t e^{(t-s)(d_2\Del-\lm_2)} u(s)ds\right\|_p
\le\,&C(p) \int_0^t e^{-\lm_2(t-s)}(t-s)^{-2\left(1-\frac1p\right)} \|u(s)\|_1 ds
\\
=\,& C(p) \|u_0\|_1\mathcal{I}_{\lm_2,p}(t).
\end{align*}
Hence the $L^p$-norm of $w$ is estimated as
\begin{align*}
\|w(t)\|_p\le\,& \| e^{t(d_2\Del-\lm_2)} w_0\|_p+\left\|\int_0^t e^{(t-s)(d_2\Del-\lm_2)} u(s)ds\right\|_p
\\
\le\,& C(p)(\|w_0\|_p+\|u_0\|_1\mathcal{I}_{\lm_2,p}(t)).
\end{align*}
Similarly for $2\le q <\infty$ we have
\begin{align*}
\|v(t)\|_q\le\,& \| e^{t(d_1\Del-\lm_1) } v_0\|_q+\int_0^t \| e^{(t-s)(d_1\Del-\lm_1)} w(s) \|_qds
\\
\le\,& C(q)\|v_0\|_q+C(q) \int_0^t e^{-\lm_1(t-s)} (t-s)^{-2\left(\frac{1+q}{2q}-\frac1q\right)}\|w(s)\|_{\frac{2q}{q+1}} ds
\\
\le\,& C(q)\left(\|v_0\|_q+\|w_0\|_{\frac{2q}{q+1}}\,\mathcal{I}_{\lm_1,q}(t)+\|u_0\|_1 \int_0^t e^{-\lm_1(t-s)} (t-s)^{-1+\frac1q}\mathcal{I}_{\lm_2,{\frac{2q}{1+q}}}(s) ds\right)
\\
\le\,& C(q) \left(\|v_0\|_q + \|w_0\|_{\frac{2q}{q+1}}\,\mathcal{I}_{\lm_1,q}(t)+\|u_0\|_1 \,\mathcal{I}_{\lm_1,q}(t)\,\mathcal{I}_{\lm_2,{\frac{2q}{1+q}}}(t) \right).
\end{align*}
Moreover it follows that
\begin{align*}
\|\N v (t)\|_2 \le\,&\| \N v_0\|_2 + C \int_0^t e^{-\lm_1(t-s)} (t-s)^{-2\left(\frac23-\frac12\right)-\frac12}\|w(s)\|_{\frac32} ds
\\
\le\,&C\| \N v_0\|_2 + \|w_0\|_{\frac32} \mathcal{I}_{\lm_1,6}(t) + \|u_0\|_1\,\mathcal{I}_{\lm_1,6}(t)\,\mathcal{I}_{\lm_2,\frac{3}{2}}(t)),
\end{align*}
which ends the proof.
\end{pr}

\vspace{5mm}
\section{Existence of  global solutions}
\vspace{5mm}
In this section, we prove the existence of a global solution to~\eqref{eqn;FS-model} as reported in Theorem~\ref{thm;GE}. To this end, introducing the chemical energy minimization and the Brezis--Merle inequality for the 4th order elliptic equation,  we show that the mass constraint~\eqref{eqn;global-assumption} leads to additional estimates for the solution to~\eqref{eqn;FS-model} and we proceed after that to the derivation of $L^p$-estimates ensuring global existence.

Throughout this section, we fix $(u_0,v_0,w_0)\in L_+^1(\R^4)\times W_+^{2,2}(\R^4)\times L_+^2(\R^4)$ such that $(1+u_0)\log(1+u_0)\in L^1(\R^4)$. We denote the corresponding solution to~\eqref{eqn;FS-model} by $(u,v,w)$ with maximal existence time $T$. Also, $C$ and $(C_i)_{i\ge 1}$ denote positive constants depending only on $(d_1,d_2,\lambda_1,\lambda_2)$ and $(u_0,v_0,w_0)$, but not on $T$, that may vary from line to line. The dependence of $C$ and $(C_i)_{i\ge 1}$ upon additional parameters will be indicated explicitly.
\subsection{Modified Lyapunov functional }

The system~\eqref{eqn;FS-model} has the Lyapunov functional $\widetilde{\F}$ defined by
\begin{align*}
\widetilde{\F}(u,v):=
\int_{\R^4}u\log u dx
+\frac{1}{2}\int_{\R^4}|\pt_t v|^2dx
+\Ep(v;u),
\end{align*}
where $\Ep(v;u)$ is the chemical energy given as
\begin{align*}
\Ep(v;u):=\frac{d_1d_2}{2}\int_{\R^4}|\Del v|^2dx
+\frac{d_1\lm_2+d_2\lm_1}{2}\int_{\R^4}|\N v|^2 dx
+\frac{\lm_1\lm_2}{2}\int_{\R^4}v^2 dx-\int_{\R^4}uvdx.
\end{align*}
Then $\widetilde{\F}(u,v)$ satisfies that
\begin{align*}
\widetilde{\F}(u,v)+\int_0^t \widetilde{\D}(u(s),v(s))ds=\widetilde{\F}(u_0,v_0),
\end{align*}
where
\begin{align*}
\widetilde{\D}(u,v):=
\int_{\R^4}u|\N(\log u - v)|^2dx+(d_2+d_1)\int_{\R^4}|\pt_t\N v|^2 dx
+(\lm_2+\lm_1)\int_{\R^4}|\pt_t v|^2dx,
\end{align*}
see Fujie--Senba \cite{Fu-Se17}. However, using the Lyapunov functional $\widetilde{\F}$ on the whole space requires at least the weight assumption $u_0 \in L^1(\R^4, \log(1+|x|^4)dx)$, so that the entropy $u\log u \in L^1(\R^4)$.
In order to get rid of the weight assumption, let us introduce the modified Lyapunov functional  $\F(u,v)$ for the problem~\eqref{eqn;FS-model},
 which plays an essential role in showing the global existence proof,
\begin{align*}
\F(u,v):=\,&
\int_{\R^4}(1+u)\log (1+u) dx
+\frac{1}{2}\int_{\R^4}|\pt_t v|^2dx
+\Ep(v;u).
\eqntag
\label{eqn;Lyapunov}
\end{align*}
The modified entropy $\int_{ \R^4}(1+u)\log(1+u)dx$ has the advantage of being nonnegative, even though it is no longer a Lyapunov functional. Indeed, $\F$ satisfies the following identity:
\vspace{5mm}
\begin{lem}\label{lem;lyapunov2}
For the solution $(u,v,w)$ to~\eqref{eqn;FS-model},
let $\F(u,v)$ be the modified Lyapunov functional  for the problem~\eqref{eqn;FS-model} defined in~\eqref{eqn;Lyapunov}.
	Then the following identity holds:
	\begin{align*}
	\F(u(t),v(t)) + \int_0^t\D(u(s),v(s))ds =\F(u_0,v_0)+\frac{1}{4}\int_0^t\int_{\R^4}|\N v(s)|^2 dxds
	\end{align*}
for $0<t<T$,	where
	\begin{align*}
	\D(u,v)=\,&\int_{\R^4}u|\N(\log (1+u) - v)|^2dx
	+\int_{\R^4}|\N(\log (1+u) -\frac12 v)|^2dx
	\\
	&+(d_2+d_1)\int_{\R^4}|\pt_t\N v|^2 dx
	+(\lm_2+\lm_1)\int_{\R^4}|\pt_t v|^2dx.
	\end{align*}
\end{lem}

\vspace{5mm}
\begin{pr}{Lemma~\ref{lem;lyapunov2}}
	Differentiating with respect to time $t$ and using the first equation in~\eqref{eqn;FS-model}, we have
	\begin{align*}
	\frac{d}{dt}\int_{\R^4}(1+u)\log(1+u)dx
	=\,&
	-\int_{\R^4}(1+u)|\N\log(1+u)|^2dx
	\\
	&+\int_{\R^4} u\N\log(1+u)\cd\N vdx.
	\end{align*}
Also, using again the first equation in~\eqref{eqn;FS-model},
	\begin{align*}
	\frac{d}{dt}\int_{\R^4}uv dx
	=\,&\int_{\R^4}v\pt_t u  dx + \int_{\R^4}u\pt_t v dx
	\equiv\,I_1+I_2,
	\end{align*}
and the first term $I_1$ is computed as
	\begin{align*}
	I_1=-\int_{\R^4}(1+u)\N\log(1+u)\cd\N vdx
	+\int_{\R^4}u|\N v|^2dx.
	\end{align*}
	Using the second and third equations of~\eqref{eqn;FS-model},
	we also have
	\begin{align*}
	I_2=\,& \int_{\R^4}\pt_t v \big(\pt_t w-d_2\Del w+\lm_2 w\big)dx
	\\
	=\,&
	\int_{\R^4}\pt_t v \pt_t w dx 
	+\int_{\R^4}\big(\pt_t v-d_1\Del v+\lm_1 v\big)\big(-d_2\Del+\lm_2\big)[\pt_t v] dx
	\\
	=\,&
	\int_{\R^4}\pt_t v \pt_t \big(\pt_t v-d_1\Del v+\lm_1 v\big) dx 
	+\int_{\R^4}\pt_t v\big(-d_2\Del+\lm_2\big)[\pt_t v] dx
	\\
	&+\int_{\R^4}\big(-d_1\Del+\lm_1\big)[v]\big(-d_2\Del+\lm_2\big)[\pt_t v] dx
	\\
	=\,&
	\frac{1}{2}\frac{d}{dt}\int_{\R^4}|\pt_t v|^2dx
	+(d_2+d_1)\int_{\R^4}|\pt_t\N v|^2dx
	\\
	&+(\lm_2+\lm_1)\int_{\R^4}|\pt_t v|^2dx
	+\frac{1}{2}\frac{d}{dt}\int_{\R^4}\big(-d_1\Del+\lm_1\big)[v]\big(-d_2\Del+\lm_2\big)[v] dx.
	\end{align*}
	Hence, we obtain
	\begin{align*}
	\frac{d}{dt}\int_{\R^4}uv dx
	=\,&
	-\int_{\R^4}(1+u)\N\log(1+u)\cd\N vdx
	+\int_{\R^4}u|\N v|^2dx
	\\
	&+\frac{1}{2}\frac{d}{dt}\int_{\R^4}|\pt_t v|^2dx
	+(d_2+d_1)\int_{\R^4}|\pt_t\N v|^2dx
	\\
	&+(\lm_2+\lm_1)\int_{\R^4}|\pt_t v|^2dx
	+\frac{1}{2}\frac{d}{dt}\int_{\R^4}\big(-d_1\Del+\lm_1\big)[v]\big(-d_2\Del+\lm_2\big)[v] dx,
	\end{align*}
so that
	\begin{align*}
	\frac{d}{dt}\left( \frac12\int_{ \R^4}|\pt_t v|^2 dx+\Ep(v;u)\right)
	=\,&\int_{\R^4}(1+u)\N\log(1+u)\cd\N vdx
	-\int_{\R^4}u|\N v|^2dx
	\\&-(d_2+d_1)\int_{\R^4}|\pt_t\N v|^2dx-(\lm_2+\lm_1)\int_{\R^4}|\pt_t v|^2dx.
	\end{align*}
Combining the above estimates, we have
\begin{align*}
&\frac{d}{dt}\left(\int_{ \R^4}(1+u)\log(1+u)dx+\frac12\int_{ \R^4}|\pt_t v|^2dx+\Ep(v;u)\right)
\\
=\,&-\int_{\R^4}(1+u)|\N\log(1+u)|^2dx+\int_{\R^4}u\N\log(1+u)\cd\N vdx
\\
&+\int_{\R^4}(1+u)\N\log(1+u)\cd\N vdx
-\int_{\R^4}u|\N v|^2dx
\\&-(d_2+d_1)\int_{\R^4}|\pt_t\N v|^2dx-(\lm_2+\lm_1)\int_{\R^4}|\pt_t v|^2dx
\\
=\,&-\int_{ \R^4}u|\N(\log(1+u)-v)|^2dx
-\int_{ \R^4}|\N(\log(1+u)-\frac12 v)|^2dx+\frac14\int_{ \R^4}|\N v|^2dx
\\
&-(d_2+d_1)\int_{\R^4}|\pt_t\N v|^2dx-(\lm_2+\lm_1)\int_{\R^4}|\pt_t v|^2dx,
\end{align*}
and the proof is complete.
\end{pr}


\vspace{5mm}
\subsection{Application of the modified Lyapunov functional}
Applying the chemical energy minimization and the Brezis--Merle type inequality leads us to the following estimates under the additional mass constraint~\eqref{eqn;global-assumption} on $u_0$.
\vspace{5mm}
\begin{prop}\label{prop;apriori-estimate}
Assume that the initial condition $u_0$ satisfies~\eqref{eqn;global-assumption}.
	Then the following estimate holds:
	\begin{align*}
	\int_{\R^4}(1+u(t))\log(1+u(t))dx+\frac{1}{2}\int_{\R^4}|\pt_t v(t)|^2 dx
	+\int_0^t \D(u(s),v(s)) ds
	\le C(\tau)
	\end{align*}
for $t\in [0,\tau]\cap [0,T)$, where $\D(u,v)$ is the functional defined in Lemma~\ref{lem;lyapunov2}.
\end{prop}
\vspace{5mm}
The functional $\Ep(v;f)$ is actually the energy
associated to the 4th order elliptic equation
$(-d_1\Del+\lm_1)(-d_2\Del+\lm_2)v=f$; that is,
\begin{align*}
\E(v;f):=\,&
\frac{d_1d_2}{2}\int_{\R^4}|\Del v|^2dx
+\frac{d_1\lm_2+d_2\lm_1}{2}\int_{\R^4}|\N v|^2 dx
+\frac{\lm_1\lm_2}{2}\int_{\R^4}v^2 dx-\int_{\R^n}fvdx.
\eqntag
\label{def;chemical-energy}
\end{align*}
The forthcoming detailed analysis of this functional, inspired from \cite[Lemma~2.2]{Ca-Co08}, plays an important role in
obtaining a global solution under the optimal mass constraint~\eqref{eqn;global-assumption}. For further use, given $f\in L^1(\R^4)$, we recall that there is a unique solution $v_f$ to 
\begin{equation*}
	\left\{
	\begin{aligned}
		&(-d_1\Del+\lm_1)(-d_2\Del+\lm_2)v_f= f,&x\in\R^4,
		\\
		&\quad\quad v_f(x) \to 0,&|x|\to\infty.
	\end{aligned}
	\right.
	\eqntag\label{eqn;BM-eq}
	\end{equation*} 
\vspace{5mm}

\begin{lem}[The chemical energy minimization]\label{lem;Chemical-EM}
Let $f\in L_+^1(\R^4)\cap  L^2(\R^4)$ and $v\in W^{2,2}(\R^4)$. Then the chemical energy $\mathcal{E}(v;f)$  defined by~\eqref{def;chemical-energy} is finite and satisfies 
	\begin{align*}
	&\E(v;f)-\E(v_f;f)
	\\
	&\qquad =\frac{d_1d_2}{2}\int_{\R^4}|\Del (v-v_f)|^2dx
	+\frac{d_1\lm_2+d_2\lm_1}{2}\int_{\R^4}|\N (v-v_f)|^2 dx\\
	&\qquad\qquad+\frac{\lm_1\lm_2}{2}\int_{\R^4}(v-v_f)^2 dx
	\ge0,
	\end{align*}
	and
	\begin{equation*}
	\Ep(v_f;f)=-\frac12\int_{ \R^4}v_ff dx.
	\end{equation*}
\end{lem}

\vspace{5mm}

\begin{pr}{Lemma~\ref{lem;Chemical-EM}}
	Let  $v\in W^{2,2}(\R^4)$. Then  
	\begin{align*}
	&\E(v;f)-\E(v_f;f)
	\\
	&\qquad =\frac{d_1d_2}{2}\int_{\R^4}|\Del (v-v_f)|^2dx
	+\frac{d_1\lm_2+d_2\lm_1}{2}\int_{\R^4}|\N (v-v_f)|^2 dx
	+\frac{\lm_1\lm_2}{2}\int_{\R^4}(v-v_f)^2 dx
	\\
	&\qquad\qquad +d_1d_2\int_{\R^4}\Del (v-v_f) \Del v_f dx
	+(d_1\lm_2+d_2\lm_1)\int_{\R^4}\N (v-v_f) \cd\N v_f dx
	\\
	&\qquad\qquad +\lm_1\lm_2 \int_{\R^4} (v-v_f) v_f dx - \int_{\R^4} v f dx+\int_{\R^4} v_f f dx.
	\end{align*}
	Integrating by parts, we see that
	\begin{align*}
	&d_1d_2\int_{\R^4}\Del (v-v_f) \Del v_f dx
	+(d_1\lm_2+d_2\lm_1)\int_{\R^4}\N (v-v_f) \cd\N v_f dx+\lm_1\lm_2\int_{\R^4} (v-v_f) v_fdx
	\\
	&\qquad =d_1d_2\int_{\R^4}(v-v_f)\Del^2 v_fdx
	-(d_1\lm_2+d_2\lm_1) \int_{\R^4}(v-v_f) \Del v_f dx
	+\lm_1\lm_2\int_{\R^4} (v-v_f) v_f dx
	\\
	&\qquad =\int_{\R^4} (v-v_f)(-d_1\Del+\lm_1)(-d_2\Del+\lm_2)[v_f] dx
	=\int_{ \R^4} f (v-v_f)dx,
	\end{align*}
	which implies that
	\begin{align*}
	&\E(v;f)-\E(v_f;f)
	\\
	&\qquad =\frac{d_1d_2}{2}\int_{\R^4}|\Del (v-v_f)|^2dx
	+\frac{d_1\lm_2+d_2\lm_1}{2}\int_{\R^4}|\N (v-v_f)|^2 dx
	\\&
	\qquad\qquad +\frac{\lm_1\lm_2}{2}\int_{\R^4}(v-v_f)^2 dx.
	\end{align*}
	This ends the proof.
\end{pr}
\vspace{5mm}

In order to obtain Proposition~\ref{prop;apriori-estimate},
let us recall the Brezis--Merle type inequality shown in Hosono--Ogawa \cite[Theorem 1.4]{Ho-Og} (cf. Brezis--Merle \cite[Theorem~1]{Br-Me}).

\begin{lem}[4-dimensional Brezis--Merle inequality]\label{lem;BM}
Consider $f\in L^1(\R^4)$ such that
\begin{equation*}
	\|f\|_1<32\pi^2d_1d_2,
\end{equation*}
and assume that $\lm_1,\lm_2>0$. Then there exists a constant $C>0$ independent of $\|f\|_1$  such that
	\begin{equation*}
	\int_{\R^4}(e^{v_f}-1)dx\le
	\frac{Ce^{2\gm\|f\|_1}}{\kp}\left(\frac{2^{\gm\|f\|_1}}{32\pi^2d_1d_2-\|f\|_1}+1\right),
	\end{equation*}
	where $\gm=1/(8\pi^2d_1d_2)$ and $\kp:=\min\{\lm_1/d_1, \lm_2/d_2\}$, recalling that $v_f$ is the solution to~\eqref{eqn;BM-eq} defined before Lemma~\ref{lem;Chemical-EM}.
\end{lem}
\vspace{5mm}
We now show how the smallness assumption on the initial mass~\eqref{eqn;global-assumption} of $u_0$ allows us to apply the Brezis--Merle type inequality stated in Lemma~\ref{lem;BM} to derive a lower bound the chemical energy $\Ep(v;u)$.
\vspace{5mm}
\begin{lem}\label{lem;Ep-bdd}
Assume that the initial condition $u_0$ satisfies~\eqref{eqn;global-assumption}. Then, for $\tau>0$, there exists a constant $C(\tau)>0$ depending on $\tau$ such that
\begin{equation*}
-\Ep(v(t);u(t))\le C(\tau), \qquad t\in [0,\tau]\cap [0,T).
\end{equation*}
\end{lem}
\vspace{5mm}
\begin{pr}{Lemma~\ref{lem;Ep-bdd}}
Since the Brezis--Merle type inequality stated in Lemma~\ref{lem;BM} is only available for positive $\lambda_1$ and $\lambda_2$, we split the analysis according to the positivity of these two parameters.

\vspace{5mm}
\noindent\textbf{Case~1: $\lm_1,\lm_2>0$.} Recalling that $v_u$ is defined in~\eqref{eqn;BM-eq} (with $f=u$), Lemma~\ref{lem;Chemical-EM} gives
\begin{equation*}
-\Ep(v,u) \le -\Ep(v_u;u) = \frac12\int_{ \R^4}u v_u dx,
\end{equation*}
and it follows from the definition~\eqref{eqn;Lyapunov} of $\F$ and Lemma \ref{lem;lyapunov2} that, for $t\in [0,\tau]\cap [0,T)$ and $\al>0$,
\begin{align*}
-\Ep(v(t);u(t))
=\,&\al\Ep(v(t);u(t)) -(1+\al)\Ep(v(t);u(t))
\\
\le\,&\al\Ep(v(t);u(t))  -(1+\al)\Ep(v_{u(t)};u(t))
\\
=\,&\al\left\{\F(u(t),v(t))-\int_{\R^4}(1+u(t))\log(1+u(t))dx-\frac{1}{2}\int_{\R^4}|\pt_t  v(t)|^2 dx\right\}
\\
& \qquad +\frac{1+\al}{2}\int_{ \R^4}u(t) v_{u(t)} dx
\\
\le\,&\al \F(u_0,v_0) 
+ \frac{\al}{4}\int_0^t \int_{\R^4}|\N v(s)|^2dxds
\\
&-\al\int_{ \R^4}(1+u(t))\log(1+u(t))dx
+\frac{1+\al}{2}\int_{ \R^4}u(t) v_{u(t)} dx.
\eqntag
\label{eqn;Ep-bdd}
\end{align*}
Now, Young's inequality, 
\begin{equation*}
	ab\le (1+a)b \le (1+a)\log(1+a) - (1+a) + e^b \le (1+a)\log(1+a)+e^b-1
\end{equation*} 
for $a,b\ge 0$ implies that
\begin{align*}
\frac{1+\al}{2}\int_{ \R^4} u(t) v_{u(t)} dx
\le \al \int_{ \R^4}(1+u(t))\log(1+u(t))dx+ \al \int_{ \R^4}\left(e^{ \frac{1+\al}{2\al} v_{u(t)}} -1\right)dx,
\end{align*}
so that we have from the inequality~\eqref{eqn;Ep-bdd} 
\begin{align*}
-\Ep(v(t);u(t))\le \al \F(u_0,v_0) 
+ \frac{\al}{4}\int_0^t \int_{\R^4}|\N v(s)|^2dxds
+\al \int_{ \R^4}\left(e^{ \frac{1+\al}{2\al} v_{u(t)}} -1\right)dx.
\eqntag 
\label{eqn;Ep-1}
\end{align*}
Owing to the assumption~\eqref{eqn;global-assumption} on $\|u_0\|_1$, we are in a position to employ the Brezis--Merle type inequality stated in Lemma~\ref{lem;BM} to obtain a bound on the last term on the right-hand side of the above inequality at the expense of an appropriate choice of $\alpha$. Indeed, since $(8\pi)^2 d_1 d_2 > \|u_0\|_1$ by~\eqref{eqn;global-assumption}, we may pick $\alpha>0$ large enough such that
\begin{equation*}
	1 < \frac{1+\alpha}{\alpha} < \frac{(8\pi)^2 d_1 d_2}{\|u_0\|_1}.
	\end{equation*}
Then, with this choice of $\alpha$, we deduce from Lemma~\ref{lem;BM} and the identity $\frac{1+\alpha}{2\alpha} v_{u(t)} = v_{\frac{1+\alpha}{2\alpha}u(t)}$ that
\begin{equation*}
\al\int_{\R^4}\left(e^{\frac{1+\al}{2\al} v_{u(t)}}-1\right)dx \le C.
\eqntag
\label{eqn;BM-bdd}
\end{equation*}
In addition, by Lemma~\ref{lem;vw-bound}, we have
 \begin{equation*}
 \frac{1}{4}\int_0^t \int_{\R^4}|\N v(s)|^2dxds \le C(\tau).
 \eqntag
 \label{eqn;additional-term}
 \end{equation*}
 Combining the estimates~\eqref{eqn;Ep-1}--\eqref{eqn;additional-term}, we obtain
 \begin{equation*}
 -\Ep(v(t);u(t)) \le C(\tau),
 \end{equation*}
which completes the proof of Lemma~\ref{lem;Ep-bdd} in that case.
 	
 \vspace{5mm}
 \noindent\textbf{Case~2: $\lm_1=\lm_2=0$.} In this case, unlike the case of $\lm_1,\lm_2>0$, the above argument does not work well due to
the alteration of  the behavior of $v_u$ as $|x|\to\infty$, which prevents the availability of a Brezis--Merle type inequality. Hence, it is necessary to modify the chemical energy minimization in Lemma~\ref{lem;Chemical-EM} and to introduce a different intermediate function. Specifically, given $f\in L^1(\R^4)$, let $\tilde{v}_f$ be the unique solution to 
\begin{equation*}
	\left\{
	\begin{aligned}
		&(-d_1\Del+1)(-d_2\Del+1) \tilde{v}_f = f,&x\in\R^4,
		\\
		&\quad\quad \tilde{v}_f(x) \to 0,&|x|\to\infty.
	\end{aligned}
	\right.
	\eqntag\label{nd3}
\end{equation*} 
 Now,  for $t\in [0,\tau]\cap [0,T)$ and $\alpha>0$,
 \begin{align}
 - \mathcal{E}(v(t);u(t)) & = \alpha \mathcal{E}(v(t);u(t)) - (1+\alpha) \mathcal{E}(v(t);u(t)) \nonumber\\
 & =  \alpha \mathcal{E}(v(t);u(t)) - (1+\alpha) \big[ \mathcal{E}(v(t);u(t))- \mathcal{E}(\tilde{v}_{u(t)};u(t)) \big] \nonumber\\
 & \qquad - (1+\alpha) \mathcal{E}(\tilde{v}_{u(t)};u(t)) . \label{nd2}
 \end{align}
 Observe that, by~\eqref{nd3}
 \begin{align}
 \mathcal{E}(v;u) - \mathcal{E}(\tilde{v}_u;u) 
 & = \frac{d_1 d_2}{2} \|\Delta(v-\tilde{v}_u)\|_2^2 
 + d_1 d_2 \int_{\R^4} \Delta \tilde{v}_u\Delta (v-\tilde{v}_u) dx 
 + \frac{d_1 d_2}{2} \|\Delta \tilde{v}_u\|_2^2 \nonumber \\
 & \qquad - \int_{\R^4} u v dx 
 - \frac{d_1 d_2}{2} \|\Delta \tilde{v}_u\|_2^2 
 + \int_{\R^4} u \tilde{v}_u dx \nonumber \\
 &  = \frac{d_1 d_2}{2} \|\Delta(v-\tilde{v}_u)\|_2^2 
 + d_1 d_2 \int_{\R^4}  (v-\tilde{v}_u) \Delta^2 \tilde{v}_u dx
  - \int_{\R^4} u (v-\tilde{v}_u) dx \nonumber \\
 & = \frac{d_1 d_2}{2} \|\Delta(v-\tilde{v}_u)\|_2^2 
 + \int_{\R^4}  (v-\tilde{v}_u) \big[ -  \tilde{v}_u
 +  (d_1+d_2) \Delta \tilde{v}_u \big] dx \nonumber \\
 & =  \frac{d_1 d_2}{2} \|\Delta(v-\tilde{v}_u)\|_2^2 
 +(d_1+d_2) \int_{\R^4} \tilde{v}_u \Delta (v-\tilde{v}_u) dx 
 - \int_{\R^4}  \tilde{v}_u (v-\tilde{v}_u) dx \nonumber \\
 & = \frac{d_1 d_2}{2} \left\| \Delta (v-\tilde{v}_u) 
 + \frac{(d_1+d_2)}{d_1 d_2} \tilde{v}_u \right\|_2^2 
 - \frac{ (d_1+d_2)^2}{2d_1 d_2} \|\tilde{v}_u\|_2^2 \nonumber\\
 & \qquad - \int_{\R^4} \tilde{v}_u (v-\tilde{v}_u) dx, \label{nd4}
 \end{align}
 and 
\begin{align}
\mathcal{E}(\tilde{v}_u;u) & = \frac{d_1d_2}{2}\int_{\R^4}  \tilde{v}_u\Delta^2 \tilde{v}_u dx -\int_{\R^4}  u \tilde{v}_u dx \nonumber \\
& = \frac{1}{2} \int_{\R^4}  \tilde{v}_u \big[  (d_1+d_2) \Delta \tilde{v}_u -  \tilde{v}_u - u \big]\ dx \nonumber \\
& = - \frac{1}{2} \int_{\R^4} u \tilde{v}_u dx - \frac{ (d_1+d_2)}{2} \|\nabla \tilde{v}_u\|_2^2 - \frac{1}{2} \|\tilde{v}_u\|_2^2. \label{nd5}
\end{align}
Gathering~\eqref{nd2}, \eqref{nd4}, and~\eqref{nd5} and Lemma \ref{lem;lyapunov2}, we find for $t\in [0,\tau]\cap [0,T)$ and $\alpha>0$
 \begin{align*}
 -\mathcal{E}(v(t);u(t)) & = \alpha \mathcal{F}(u(t),v(t))
  - \alpha\int_{\R^4}  (1+u(t)) \log{(1+u(t))} dx 
  - \frac{\alpha}{2} \|\partial_t v(t)\|_2^2 \\
 & \qquad - (1+\alpha) \frac{d_1 d_2}{2} \left\| \Delta \big(v(t)-\tilde{v}_{u(t)}\big) 
 + \frac{(d_1+d_2)}{d_1 d_2} \tilde{v}_{u(t)} \right\|_2^2 \\
 & \qquad + \frac{(1+\alpha)  (d_1+d_2)^2}{2d_1 d_2} \|\tilde{v}_{u(t)}\|_2^2 + (1+\alpha) \int_{\R^4}  \tilde{v}_{u(t)} \big(v-\tilde{v}_{u(t)}\big) dx \\
 & \qquad + \frac{1+\alpha}{2}\int_{\R^4}  u(t) \tilde{v}_{u(t)} dx 
 + \frac{(1+\alpha)  (d_1+d_2)}{2} \|\nabla \tilde{v}_{u(t)}\|_2^2 \\ 
 & \qquad + \frac{(1+\alpha)}{2} \|\tilde{v}_{u(t)}\|_2^2 \\
 & \le \alpha \mathcal{F}(u_0,v_0) + \frac{\alpha}{4} \int_0^t \|\nabla v(s)\|_2^2 ds 
 - \alpha \int_{\R^4}  (1+u(t)) \log{(1+u(t))} dx \\
 & \qquad + \frac{1+\alpha}{2} \int_{\R^4}  u(t) \tilde{v}_{u(t)} dx 
 + (1+\alpha)  \int_{\R^4}  v(t) \tilde{v}_{u(t)} dx \\
 & \qquad + \frac{(1+\alpha) (d_1+d_2)}{2} \|\nabla \tilde{v}_{u(t)}\|_2^2 
 + \frac{(1+\alpha)  \big( d_1^2 + d_2^2 \big)}{2 d_1 d_2} \|\tilde{v}_{u(t)}\|_2^2.
 \end{align*}
 At this point, we use Young's inequality as in the previous case to estimate
 \begin{equation*}
 \frac{1+\alpha}{2} \int_{\R^4}  u(t) \tilde{v}_{u(t)} dx \le \alpha \int_{\R^4}  (1+u(t)) \log{(1+u(t))} dx + \alpha \int_{\R^4}  \left( e^{\frac{1+\alpha}{2\alpha} \tilde{v}_{u(t)}} - 1 \right) dx
 \end{equation*}
 and thereby obtain
 \begin{align*}
 -\mathcal{E}(v(t);u(t)) & \le \alpha \mathcal{F}(u_0,v_0) 
 + \frac{\alpha}{4} \int_0^t \|\nabla v(s)\|_2^2 ds 
 + \alpha \int_{\R^4} \left( e^{\frac{1+\alpha}{2\alpha} \tilde{v}_{u(t)}} - 1 \right)\ dx \nonumber \\
 & \qquad + (1+\alpha) \int_{\R^4}  v(t) \tilde{v}_{u(t)}\ dx + \frac{(1+\alpha)  (d_1+d_2)}{2} \|\nabla \tilde{v}_{u(t)}\|_2^2 \nonumber \\
 & \qquad + \frac{(1+\alpha) \big( d_1^2 + d_2^2 \big)}{2 d_1 d_2} \|\tilde{v}_{u(t)}\|_2^2. 
 \end{align*}
 Noting that elliptic regularity and~\eqref{eqn;mass} ensure that
 \begin{align*} 
 \| \tilde{v}_{u(t)}\|_2+ \|\N \tilde{v}_{u(t)}\|_2 \le C \|u(t)\|_1 = C \|u_0\|_1,
 \end{align*}
we deduce from Lemma~\ref{lem;vw-bound} and the Brezis--Merle type inequality in Lemma~\ref{lem;BM} that, for an appropriate choice of $\alpha$ large enough,
\begin{align*}
 -\mathcal{E}(v(t);u(t))  \le C(\tau),
\end{align*}
 and this concludes the proof.
 
 \vspace{5mm}
 \noindent\textbf{Case~3: $\lambda_1>\lambda_2=0$ or $\lambda_2>\lambda_1=0$.} These two cases can be dealt with as the previous case and we thus omit the proof.
\end{pr}

\vspace{5mm}
\begin{pr}{Proposition~\ref{prop;apriori-estimate}}
By virtue of Lemma~\ref{lem;lyapunov2}, we infer from~\eqref{eqn;additional-term}, Lemma~\ref{lem;vw-bound}, and Lemma~\ref{lem;Ep-bdd} that, for $t\in [0,\tau]\cap [0,T)$,
	\begin{align*}
	&\int_{\R^4}(1+u(t))\log (1+u(t)) dx
	+\frac{1}{2}\int_{\R^4}|\pt_t v|^2dx
	+\int_0^t \D(u(s),v(s)) ds
	\\
	&=\F(u(t),v(t)) -\Ep(v(t);u(t)) + \int_0^t \D(u(s),v(s)) ds
	\\
	&\le -\Ep(v(t);u(t)) 
	+\F(u_0,v_0)+\frac{1}{4}\int_0^t \int_{\R^4}|\N v(s)|^2 dxds
	\\
	&\le C(\tau),
	\end{align*}
and the proof is complete.
\end{pr}

\vspace{5mm}
\subsection{Energy estimates}
We proceed to show $L^p$-estimates for the first component $u$ of the solution to~\eqref{eqn;FS-model} via the estimates derived in Proposition~\ref{prop;apriori-estimate}.
The first step is to derive estimates for $w$ in $ L^{\infty}\big(0,t; L^2(\R^4)\big)\cap L^2\big(0,t; W^{1,2}(\R^4)\big)$. To this end, we first report the following functional inequality which is inspired from \cite[Equation~(22)]{Bi-He-Na94}, see also \cite[Lemma~3.5]{Na-Se-Ya1}.

\vspace{5mm}
\begin{lem}\label{lem;L3/2-estimate}
Let $1\le p\le2$. There exists a positive constant $C(p)$ such that, for all $f\in L_+^1(\R^4)\cap W^{1,2}(\R^4)$ and $N\ge 1$,
\begin{align*}
\|f\|_p^p\le
\frac{C(p)}{ \left(\log(1+N^2)\right)^{2-p}}
\|(1+f)\log(1+f)\|_1^{2-p}
\left(\int_{ \R^4}\frac{|\N f|^2}{1+f}dx\right)^{2(p-1)}
+ C(p) N^{p-1}\|f\|_1.
\end{align*}
In particular, letting $p=3/2$
\begin{align*}
\|f\|_{\frac32}^{\frac32} \le 
\frac{C}{\sqrt{\log(1+N^2)}}\|(1+f)\log(1+f)\|_1^{\frac12}\int_{\R^4}\frac{|\N f|^2}{1+f}dx
+C N^{\frac12}\|f\|_1
\end{align*}
for any $f\in L_+^1(\R^4)\cap W^{1,2}(\R^4)$ and $N \ge 1$.
\end{lem}
\vspace{5mm}

\begin{pr}{Lemma~\ref{lem;L3/2-estimate}}
For $N\ge 1$, define $\vp_N\in C([0,\infty))$ by
\begin{align*}
\vp_N(r):=
\left\{
\begin{aligned}
&0,&\text{if}\quad&0\le r\le N,
\\
&2(r-N),&\text{if}\quad&N\le r \le 2N,
\\
&r, &\text{if}\quad&2N\le r <\infty.
\end{aligned}
\right.
\end{align*} 
We note that $0\le \vp_N(r)\le r$, and
\begin{align*}
r-\vp_N(r)=
\left\{
\begin{aligned}
&r,&\text{if}\quad&0\le r\le N,
\\
&2N-r,&\text{if}\quad&N\le r \le 2N,
\\
&0, &\text{if}\quad&2N\le r <\infty
\end{aligned}
\right.
\end{align*} 
with $0 \le r-\vp_N(r)\le r$ for $r\ge 0$.
Let $1\le p \le 2$ and $f \in L_+^1(\R^4)\cap W^{1,2}(\R^4)$.
On the one hand, by the Gagliardo--Nirenberg inequality
\begin{align*}
\| f\|_{2p} \le C(p) \|f\|_2^{\frac{2-p}{p}} \|\N f\|_2^{\frac{2(p-1)}{p}},
\end{align*}
so that
\begin{align*}
\| \vp_N(\sqrt{f}) \|_{2p} \le C(p) \|\vp_N(\sqrt{f})\|_2^{\frac{2-p}{p}} \|\N\vp_N(\sqrt{f})\|_2^{\frac{2(p-1)}{p}}.
\end{align*}
Since 
\begin{align*}
\|\N \vp_N(\sqrt{f})\|_2^2
=\,&\| \vp_N'(\sqrt{f}) \N \sqrt{f} \|_2^2
\\
=\,&\int_{\{N^2\le f \le 4N^2\}}\frac{|\N f|^2}{f} dx
\le\,2\int_{\R^4}\frac{|\N f|^2}{1+f}dx,
\end{align*}
and
\begin{align*}
\|\vp_N(\sqrt{f})\|_2^2
=\,&\int_{ \{f \ge N^2\}}|\vp_N(\sqrt{f}) |^2dx
\\
\le\,&\int_{ \{f \ge N^2\}} f dx\le \frac{1}{\log(1+N^2)}\int_{ \R^4}(1+f)\log(1+f)dx,
\end{align*}
we obtain
\begin{align*}
\| \vp_N(\sqrt{f}) \|_{2p}^{2p} \le 
\frac{C(p)}{ \left(\log(1+N^2)\right)^{2-p}}
\|(1+f)\log(1+f)\|_1^{2-p}
\left(\int_{ \R^4}\frac{|\N f|^2}{1+f}dx\right)^{2(p-1)}.
\eqntag
\label{eqn;L3/2estimate2}
\end{align*}
On the other hand,
\begin{align*}
\| \vp_N(\sqrt{f})-\sqrt{f} \|_{2p}^{2p} 
\le\,&  \int_{ \{f \le 2N\}} f^p dx
\le\, (2N)^{p-1}\|f\|_1.
\end{align*}
Hence, combining the outcome of the above computation with~\eqref{eqn;L3/2estimate2}, 
we find 
\begin{align*}
\|\sqrt{f}\|_{2p}^{2p}\le
\frac{C(p)}{ \left(\log(1+N^2)\right)^{2-p}}
\|(1+f)\log(1+f)\|_1^{2-p}
\left(\int_{ \R^4}\frac{|\N f|^2}{1+f}dx\right)^{2(p-1)}
+ C(p) N^{p-1}\|f\|_1.
\end{align*}
Equivalently,
\begin{align*}
\|f\|_p^p\le
\frac{C(p)}{ \left(\log(1+N^2)\right)^{2-p}}
\|(1+f)\log(1+f)\|_1^{2-p}
\left(\int_{ \R^4}\frac{|\N f|^2}{1+f}dx\right)^{2(p-1)}
+ C(p) N^{p-1}\|f\|_1,
\end{align*}
which ends the proof of Lemma~\ref{lem;L3/2-estimate}.
\end{pr}
\vspace{5mm}

\begin{lem}\label{lem;u-1/2-log}
Assume that the initial condition $u_0$ satisfies~\eqref{eqn;global-assumption}. Then, for $\tau>0$, there exists a constant $C(\tau)>0$ depending on $\tau$ such that
\begin{equation*}
\|w(t)\|_2^2+\int_0^t \|\N w(s)\|_2^2d\t+
\int_0^t \int_{ \R^4} \frac{|\N u(s)|^2}{1+u(s)} dx ds\le C(\tau), \qquad t\in [0,\tau]\cap [0,T).
\end{equation*}
\end{lem}

\vspace{5mm}
\begin{pr}{Lemma~\ref{lem;u-1/2-log}}
Let $\tau>0$. Differentiating with respect to time and integrating by parts, we infer from~\eqref{eqn;FS-model} that
\begin{align*}
&\frac{d}{dt}\int_{ \R^4}(1+u)\log(1+u)dx
+\frac{d}{dt}\|w\|_2^2
\\
=\,&
-2d_2\|\N w\|_2^2
-\int_{ \R^4}\frac{|\N u|^2}{1+u}dx
 +\int_{ \R^4}\N u\cd\N v dx
\\
&- \int_{ \R^4}\N \log(1+u)\cd\N v dx
-2\lm_2\int_{ \R^4}w^2dx+2\int_{ \R^4}uwdx        
\\
=\,&-2d_2\|\N w\|_2^2
-\int_{ \R^4}\frac{|\N u|^2}{1+u}dx
+\frac{1}{d_1}\int_{ \R^4}u(w-\lm_1v-\pt_t v)dx
\\
&-\frac{1}{d_1}\int_{ \R^4}\log(1+u)(w-\lm_1 v-\pt_t v)dx
-2\lm_2\|w\|_2^2 +2\int_{ \R^4}uwdx.                                                                                                                                                                                                                                                                                                                                                                                                                                                                                                                                                                                                                                                                                                                                                     
\end{align*}
Using the nonnegativity of $u$, $v$, $w$, we further obtain
\begin{align*}
&\frac{d}{dt}\int_{ \R^4}(1+u)\log(1+u)dx
+\frac{d}{dt}\|w\|_2^2\\
\le\,&-2d_2\|\N w\|_2^2
-\int_{ \R^4}\frac{|\N u|^2}{1+u}dx
+\left(2+\frac{1}{d_1}\right)\int_{ \R^4}uwdx+\frac{1}{d_1}\int_{ \R^4}u|\pt_t v| dx
\\
&+\frac{\lm_1}{d_1}\int_{ \R^4}v\log(1+u)dx
+\frac{1}{d_1}\int_{ \R^4}|\pt_t v| \log(1+u)dx.      
\eqntag\label{eqn;ulog-w}   
\end{align*}
Since $\log(1+r)\le \sqrt{r}$ for $r>0$, it follows from~\eqref{eqn;mass} and Proposition~\ref{prop;apriori-estimate} that the fifth and sixth terms in the right-hand side of~\eqref{eqn;ulog-w} can be estimated on $[0,\tau]\cap [0,T)$ as
\begin{align*}
\frac{\lm_1}{d_1}\int_{ \R^4}v\log(1+u)dx
+\frac{1}{d_1}\int_{ \R^4}|\pt_t v| \log(1+u)dx
\le\,&
\frac{\lm_1}{d_1}\int_{ \R^4} v u^{\frac12}dx 
+\frac{1}{d_1}\int_{ \R^4}|\pt_t v| u^{\frac12} dx
\\
\le\,&C(\|u_0\|_1+\|\pt_t v\|_2^2+\|v\|_2^2)
\\
\le\,& C(\tau).
\end{align*}
We are left with estimating the $L^1$-norm of $uw$ and $u\pt_t v$. To this end, recalling the Sobolev inequality $\|w\|_4 \le C \|\nabla w\|_2$, we infer from~\eqref{eqn;mass} and H\"{o}lder's and Young's inequalities that
\begin{align*}
\left(2+\frac{1}{d_1}\right)\int_{ \R^4}uw dx \le\,& C\|u\|_{\frac43}\|w\|_4 \le\, C \|u\|_1^{\frac14} \|u\|_{\frac32}^{\frac34} \|\N w\|_2 \\
\le\,& C\|u\|_{\frac32}^{\frac32} + d_2  \|\N w\|_2^2.
\end{align*}
Also, according to the Gagliardo--Nirenberg inequality and Proposition~\ref{prop;apriori-estimate},
\begin{equation*}
	\|\partial_t v\|_3 \le C \|\partial_t v\|_2^{\frac13} \|\N\partial_t v\|_2^{\frac{2}{3}} \le C(\tau) \|\N\partial_t v\|_2^{\frac{2}{3}}
\end{equation*}
on $[0,\tau]\cap [0,T)$, so that, by H\"{o}lder's and Young's inequalities, 
\begin{align*}
\frac{1}{d_1}\int_{ \R^4}u|\pt_t v |dx
\le\,& \|u\|_{\frac32} \|\pt_t v\|_{3} \le C(\tau) \|u\|_{\frac32} \|\N\partial_t v\|_2^{\frac{2}{3}}
\\
\le\,& \|u\|_{\frac32}^{\frac32} + C(\tau)\|\N \pt_t v\|_2^2.
\end{align*}
Collecting the above estimates and using Lemma~\ref{lem;L3/2-estimate}, it follows from~\eqref{eqn;ulog-w} that, on $[0,\tau]\cap [0,T)$,
\begin{align*}
&\frac{d}{dt}\int_{ \R^4}(1+u)\log(1+u)dx
+\frac{d}{dt}\|w\|_2^2
+2d_2\|\N w\|_2^2
+\int_{ \R^4}\frac{|\N u|^2}{1+u}dx
\\
\le\,&
C \|u\|_{\frac32}^{\frac32}
+d_2  \|\N w\|_2^2
+ C(\tau) \big( 1 + \|\N \pt_t v\|_2^2\big)
\\
\le\,&\frac{C}{\sqrt{\log(1+N^2)}} \|(1+u)\log(1+u)\|_1^{\frac12}\int_{\R^4}\frac{|\N u|^2}{1+u}dx
+C N^{\frac12}\|u_0\|_1
\\
&+d_2  \|\N w\|_2^2 + C(\tau) \big( 1 + \|\N \pt_t v\|_2^2\big).
\end{align*}
Hence, by Proposition~\ref{prop;apriori-estimate}, 
\begin{align*}
	&\frac{d}{dt}\int_{ \R^4}(1+u)\log(1+u)dx
	+\frac{d}{dt}\|w\|_2^2
	+d_2\|\N w\|_2^2
	+\int_{ \R^4}\frac{|\N u|^2}{1+u}dx
	\\
	\le\,&\frac{C_1(\tau)}{\sqrt{\log(1+N^2)}} \int_{\R^4}\frac{|\N u|^2}{1+u}dx
		+C N^{\frac12} + C(\tau) \big( 1 + \|\N \pt_t v\|_2^2\big)
\end{align*}
on $[0,\tau]\cap [0,T)$ and taking $N(\tau)$ large enough such that
\begin{align*}
\frac{C_1(\tau)}{\sqrt{\log(1+N(\tau)^2)}} \le \frac12,
\end{align*}
it follows from Proposition~\ref{prop;apriori-estimate} by integrating with respect to time that, for $t\in [0,\tau]\cap [0,T)$,
\begin{align*}
&\int_{ \R^4}(1+u(t))\log(1+u(t))dx
+\|w(t)\|_2^2
+d_2\int_0^t\|\N w(s)\|_2^2ds
+\frac12\int_0^t\int_{ \R^4}\frac{|\N u(s)|^2}{1+u(s)}dxds
\\
\le\,&C(\tau) \int_0^t \big( 1 + \|\N \pt_t v(s)\|_2^2 \big)ds \le C(\tau),
\end{align*}
which ends the proof.
\end{pr}

\vspace{5mm}
The $L^p$-estimates of $u$ then works well by use of the standard energy method and the Nash--Moser iteration technique.

\vspace{5mm}
\begin{pr}{Theorem~\ref{thm;GE}}
Let $\tau>0$. According to Proposition~\ref{prop;apriori-estimate}, Lemma~\ref{lem;u-1/2-log}, and the definition of $\D$,
\begin{equation}
	\|\pt_t v(t)\|_2^2  +  \|w(t)\|_2^2 + \int_0^t \big( \|\N\pt_t v(s)\|_2^2 + \|\N w(s)\|_2^2 \big)ds \le C(\tau) \label{nd10}
\end{equation}
for $t\in [0,\tau]\cap [0,T)$. We now estimate the $L^p$-norm of $u$ for $p\in (1,\infty)$. Multiplying the first equation of~\eqref{eqn;FS-model} by $pu^{p-1}$ for $p>\frac43$ and integrating by parts, we obtain
\begin{align*}
\frac{d}{dt}\|u(t)\|_p^p
+\frac{4(p-1)}{p}\|\N u^{\frac p2}\|_2^2
=\,&
p(p-1)\int_{\R^n}u^{p-1}\N u \cd \N vdx
\\
=\,&
\frac{p-1}{d_1}\int_{\R^n}u^p(w-\pt_t v -\lm_1v)dx
\\
\le\,&
\frac{p-1}{d_1}\int_{\R^n}u^pwdx
+\frac{p-1}{d_1}\int_{\R^n}u^p|\pt_t v| dx,
\end{align*}
recalling that $u^pv\ge 0$. Since the H\"{o}lder and Sobolev inequalities imply that
\begin{align*}
\|f\|_{\frac83}\le \|f\|_{2}^{\frac12}\|f\|_{4}^{\frac12}
\le C\|f\|_{2}^{\frac12}\|\N f\|_{2}^{\frac12},
\end{align*}
another use of H\"older's inequality, along with that of the Sobolev and Young inequalities, gives, for any $\ep>0$,
\begin{align*}
\int_{\R^4}u^pw dx
\le\,& 
\| u^{\frac p2} \|_{\frac83}^2\|w\|_4 \le C \|u^{\frac p2}\|_{2}\|\N u^{\frac p2}\|_2 \|\N w\|_{2}
\\
\le\,&\ep\|\N u^{\frac p2}\|_2^2+C(\ep)\|\N w\|_{2}^2\|u\|_p^p.
\end{align*}
Similarly,
\begin{align*}
\int_{\R^4}u^p|\pt_t v| dx\le\ep\|\N u^{\frac p2}\|_2^2+C(\ep)\|\N \pt_t v\|_{2}^2\|u\|_p^p.
\end{align*}
Hence, 
\begin{equation*}
\begin{aligned}
\frac{d}{dt}\|u(t)\|_p^p + &\, \left(\frac{4(p-1)}{p}-\frac{2\ep(p-1)}{d_1}\right)\|\N u^{\frac p2}\|_2^2 \\
\le\,& 
C(\ep) (p-1)\|u\|_p^p\left(\|\N w\|_2^2+\|\N \pt_t  v\|_{2}^2\right),
\end{aligned}
\eqntag
\label{eqn;L^p-u}
\end{equation*}
and choosing $\ep=d_1/p$,
\begin{align*}
\frac{d}{dt}\|u(t)\|_p^p+\frac{2(p-1)}{p}\|\N u^{\frac p2}\|_2^2
\le\,& 
C(p)\|u\|_p^p\left(\|\N w\|_2^2+\|\N \pt_t  v\|_{2}^2\right).
\end{align*}
At this point, we fix $t_0\in (0,T)$ and assume further that $\tau>t_0$. Since $u(t_0)\in L^p(\R^4)$ by Proposition~\ref{prop;LWP},
 the Gronwall inequality yields by~\eqref{eqn;L^p-u} that
\begin{align*}
\|u(t)\|_p^p
\le \|u(t_0)\|_p^p\exp\left(C(p) \int_{t_0}^{t}\left(\|\N w(s)\|_2^2+\|\N \pt_t  v(s)\|_{2}^2\right)ds\right)
\end{align*}
for $t\in [t_0,\tau]\cap [t_0,T)$. It then follows from~\eqref{nd10} and the above inequality that
\begin{equation*}
\|u(t)\|_p\le C(p,t_0,\tau), \qquad t\in [t_0,\tau]\cap [t_0,T).
\eqntag
\label{eqn;Lp-u}
\end{equation*}
Parabolic regularity theory then allows us to spread the established integrability property~\eqref{eqn;Lp-u} of $u$ to $w$ and $v$. Indeed, we infer from~\eqref{eqn;FS-model} and~\eqref{eqn;Lp-u} that, for $t\in [t_0,\tau]\cap [t_0,T)$ and $p\in [1,\infty)$,
\begin{align*}
\|w(t)\|_p\le\,& \|e^{(t-t_0)(d_2\Del-\lm_2)} w(t_0)\|_p
+\int_{t_0}^t \|e^{(t-s)(d_2\Del-\lm_2)} u(s)\|_pds
\\
\le\,&C(p,t_0) \|w(t_0)\|_p + C(p,t_0) \int_{t_0}^t \|u(s)\|_{p} ds\le C(p,t_0,\tau),
\eqntag
\label{eqn;Lp-w}
\end{align*}
while
\begin{align*}
	\|w(t)\|_\infty\le\,& \|e^{(t-t_0)(d_2\Del-\lm_2)} w(t_0)\|_\infty
	+\int_{t_0}^t \|e^{(t-s)(d_2\Del-\lm_2)} u(s)\|_\infty ds
	\\
	\le\,&C(t_0) \|w(t_0)\|_\infty + C(t_0) \int_{t_0}^t (t-s)^{-\frac12}\|u(s)\|_{4} ds\le C(t_0,\tau).
	\eqntag
	\label{eqn;Linfty-w}
\end{align*}
Furthermore, thanks to~\eqref{eqn;Lp-w} and~\eqref{eqn;Linfty-w}, we find for $t\in [t_0,\tau]\cap [t_0,T)$ and $p\in [2,\infty]$ that
\begin{equation*}
\|\N v(t)\|_{2p}
\le\, C(p) \|\N v(t_0)\|_{2p} 
+ C(p)\int_{t_0}^{t} (t-s)^{-2\left(\frac{1}{p}-\frac{1}{2p}\right)-\frac12} \|w(s)\|_p ds
\le\, C(p,t_0,\tau).
\end{equation*}
Consequently, for all $p\in [1,\infty)$
\begin{equation*}
	\|(u\N v)(t)\|_p \le C(p,t_0,\tau), \qquad t\in [t_0,\tau]\cap [t_0,T),
\end{equation*}
and choosing $p>6$, we proceed as in Tao--Winkler \cite[Lemma~A.1]{Ta-Wi12} and perform a Moser iteration to conclude that
\begin{equation*}
\|u(t)\|_{\infty} \le C(t_0,\tau), \qquad t\in [t_0,\tau]\cap [t_0,T).
\end{equation*}
In particular, this bound excludes the finite time blow-up of the solution $(u,v,w)$ to~\eqref{eqn;FS-model}, and hence the proof is complete.
\end{pr}

\vspace{5mm}
\section{Boundedness of solutions}
This section deals with the proof of the boundedness of solutions
to~\eqref{eqn;FS-model} stated in Theorem~\ref{thm;Bdd}. We thus consider initial conditions 
\begin{equation*}
	(u_0,v_0,w_0)\in L_+^1(\R^4)\times L_+^1(\R^4)\cap W^{2,2}(\R^4) \times L_+^1(\R^4)\cap L^2(\R^4)
\end{equation*} 
with $(1+u_0)\log{(1+u_0)}\in L^1(\R^4)$ satisfying~\eqref{eqn;bdd-assumption}; that is,
\begin{equation*}
	\|u_0\|_1<\frac1{\sqrt{3}}(8\pi)^2d_1d_2 < (8\pi)^2d_1d_2,
\end{equation*}
 and recall that, since~\eqref{eqn;bdd-assumption} implies~\eqref{eqn;global-assumption}, the corresponding solution $(u,v,w)$ to~\eqref{eqn;FS-model} is global according to Theorem~\ref{thm;GE}. Throughout this section, $C$ and $(C_i)_{i\ge 1}$ denote positive constants depending only on $(d_1,d_2,\lambda_1,\lambda_2)$ and $(u_0,v_0,w_0)$, but not on time, that may vary from line to line. The dependence of $C$ and $(C_i)_{i\ge 1}$ upon additional parameters will be indicated explicitly.

\vspace{5mm}
\subsection{Time-independent energy estimates}
Let us introduce the following energy functional $\L$ for the problem~\eqref{eqn;FS-model}, which plays an essential role in ensuring uniform bounded estimates:
\begin{align*}
\L(u,v,w):=\L_0(u,v)+\frac{1}{2d_1}\|w\|_2^2+\frac{\lm_2}{d_1d_2}\| \N W \|_2^2,
\eqntag
\label{def;L_0}
\end{align*}
where 
\begin{align*}
\L_0(u,v):=\int_{ \R^4}(1+u)\log(1+u)dx
+\frac{1}{2d_1}\Ep_0(v),
\eqntag
\label{eqn;L_0}
\end{align*}
\begin{align*}
\Ep_0(v):=\|\pt_t v\|_2^2+d_1d_2\|\Del v\|_2^2+(d_1\lm_2+d_2\lm_1)\|\N v\|_2^2+\lm_1\lm_2\|v\|_2^2,
\end{align*}
and $W:=(-\Del)^{-1}w=E_4*w$ with the Poisson kernel in $\R^4$ given by $E_4(x)=1/(4\pi^2|x|^2)$ for $x\in\R^4$. In order to establish such estimates,
we derive a differential inequality for the energy functional $\L$ as follows.
\vspace{5mm}
\begin{prop}\label{prop;bound-L}
There exist constants $C_1, C_2>0$ such that
\begin{equation*}
	\frac{d}{dt}\L(u,v,w)+ C_1 \L(u,v,w) + C_1 \D_1(u,v,w)\le C_2, \qquad t\ge 0,
\end{equation*}
where $\L(u,v,w)$ is the energy functional defined in~\eqref{def;L_0},
\begin{align*}
	\D_1(u,v,w):=\D_0(u,v,w)+\frac{1}{d_1d_2}\|\N\pt_t W\|_2^2
	+\frac{\lm_2}{d_1}\|w\|_2^2
	+\frac{\lm_2^2}{d_1d_2}\|\N W\|_2^2
	\eqntag
	\label{eqn;dissipative}
\end{align*}
and
\begin{align*}
	\D_0(u,v):=\int_{ \R^4}\frac{|\N u|^2}{1+u} dx+\frac{\lm_2+\lm_1}{d_1}\|\pt_t v\|_2^2
		+\frac{d_2+ d_1}{d_1}\|\N \pt_t v\|_2^2
\end{align*} 
with $W:=E_4*w$.
\end{prop}

\vspace{5mm}
We begin with a differential inequality for the functional $\L_0$ defined in~\eqref{eqn;L_0}.
\vspace{5mm}
\begin{lem}\label{lem;L_0}
The functional $\L_0$ defined in~\eqref{eqn;L_0} satisfies the differential inequality
\begin{align*}
\frac{d}{dt}\L_0(u,v)+\D_0(u,v) \le \frac{1}{d_1}\int_{ \R^4} uw dx 
+ \frac{\lm_1}{d_1}\|u_0\|_1^{\frac12}\|v\|_2
+ \frac{1}{d_1}\|u_0\|_1^{\frac12}\|\pt_t v\|_2.
\end{align*}
\end{lem}
\vspace{5mm}
\begin{pr}{Lemma~\ref{lem;L_0}}
Applying the first equation of~\eqref{eqn;FS-model} and integrating by parts, we have
\begin{align*}
\frac{d}{dt}\int_{\R^4}(1+u)\log(1+u)dx
=\,&-\int_{\R^4}\frac{|\N u|^2}{1+u}dx+\int_{\R^4}(\log(1+u)-u)\Del vdx.
\eqntag
\label{eqn;ulogu-energy}
\end{align*}
At this point, we note from the equation for $v$ that
\begin{align*}
-\int_{\R^4}u\Del v dx=\,&-\frac{1}d_1\int_{\R^4}u(\pt_t v+\lm_1 v-w)dx
\\
=\,&-\frac{1}{d_1}\int_{\R^4}u\pt_t vdx-\frac{\lm_1}{d_1}\int_{\R^4}uvdx+\frac{1}{d_1}\int_{\R^4}uwdx,
\end{align*}
and by use of the equation of $w$ it follows that
\begin{align*}
-\frac{1}{d_1}\int_{\R^4}u\pt_t vdx
=\,&
-\frac{1}{d_1}\int_{\R^4}(\pt_tw-d_2\Del w+\lm_2 w) \pt_t vdx
\\
=\,&-\frac{1}{d_1}\int_{\R^4}(\pt_t-d_2\Del+\lm_2)(\pt_t-d_1\Del+\lm_1)[v]\, \pt_tv dx
\\
=\,&
-\frac{1}{d_1}
\Bigg[
\frac{1}{2}\frac{d}{dt}\|\pt_t v\|_2^2
+\frac{d_1d_2}{2}\frac{d}{dt}\|\Del v\|_2^2
+\frac{d_1\lm_2+d_2\lm_1}{2}\frac{d}{dt}\|\N v\|_2^2
\\
&+\frac{\lm_1\lm_2}{2}\frac{d}{dt}\|v\|_2^2
+(\lm_2+\lm_1)\|\pt_t v\|_2^2
+(d_2+d_1)\|\N\pt_t v \|_2^2
\Bigg],
\end{align*}
which shows that
\begin{align*}
-\int_{\R^4}u\Del v dx=\,&
-\frac{1}{2d_1}\frac{d}{dt}\Ep_0(v)
-\frac{1}{d_1}\Big[(\lm_2+\lm_1)\|\pt_t v\|_2^2
+(d_2+d_1)\|\N\pt_t v\|_2^2\Big]
\\
&-\frac{\lm_1}{d_1}\int_{\R^4}uvdx+\frac{1}{d_1}\int_{\R^4}uwdx.
\eqntag
\label{eqn;ulogu-energy2}
\end{align*}
Moreover, since $\log(1+r)\le \sqrt{r}$ for $r\ge0$, we deduce from~\eqref{eqn;mass}, H\"{o}lder's inequality, and the nonnegativity of $u$ and $w$ that
\begin{align*}
\int_{\R^4}\log(1+u)\Del v dx
=\,&
\frac{1}{d_1}\int_{\R^4}\log(1+u)(\pt_t v+\lm_1v-w)dx
\\
\le\,&
\frac{1}{d_1}\int_{ \R^4}u^\frac12|\pt_t v| dx
+\frac{\lm_1}{d_1}\int_{\R^4}u^{\frac12}v dx
\\
\le\,&
\frac{1}{d_1}\|u_0\|_1^{\frac12}\|\pt_t v\|_2
+\frac{\lm_1}{d_1}\|u_0\|_1^{\frac12}\|v\|_2.
\eqntag
\label{eqn;ulogu-energy3}
\end{align*}
Thus we infer from~\eqref{eqn;ulogu-energy}--\eqref{eqn;ulogu-energy3} and the nonnegativity of $u$ and $v$ that 
\begin{align*}
\frac{d}{dt}\L_0(u,v)+\D_0(u,v)
\le\,&
\frac{1}{d_1}\int_{\R^4}uw dx
+\frac{1}{d_1}\|u_0\|_1^{\frac12}\|\pt_t v\|_2
+\frac{\lm_1}{d_1}\|u_0\|_1^{\frac12}\|v\|_2.
\end{align*}
The proof is complete.
\end{pr}

\vspace{5mm}
We now estimate the $L^1$-norm of $uw$,
and here utilize the  Hardy--Littlewood--Sobolev inequality~\eqref{eqn;HLS-inequality} and its optimal constant $C_{HLS} =(3/2)^{\frac12}\pi$.
\begin{lem}\label{lem;estimate-uw}
There holds
\begin{align*}
\frac{1}{d_1}\int_{\R^4}uw dx\le\,&
\frac{1}{4\pi^2d_1d_2}C_{HLS}\|u\|_{\frac43}^2
-\frac{1}{2d_1}\frac{d}{dt}\|w\|_2^2
-\frac{\lm_2}{d_1d_2}\frac{d}{dt}\|\N W\|_2^2
\\
&-\frac{1}{d_1d_2}\|\N\pt_t W\|_2^2
-\frac{\lm_2}{d_1}\|w\|_2^2
-\frac{\lm_2^2}{d_1d_2}\|\N W\|_2^2,
\end{align*}
with $C_{HLS}=(3/2)^{\frac12}\pi$ and $W=(-\Del)^{-1}w=E_4*w$.
\end{lem}
\vspace{5mm}

\begin{pr}{Lemma~\ref{lem;estimate-uw}}
We set $U:=(-\Del)^{-1} u = E_4*u$, so that $-\Del U=u$ in $\R^4$.
 Integrating by parts, we observe that
\begin{align*}
\frac{1}{d_1}\int_{\R^4}uw dx
=\,&
-\frac{1}{d_1}\int_{\R^4}w\Del U dx
\\
=\,&-\frac{1}{d_1d_2}\int_{\R^4}U(\pt_t w+\lm_2 w-u)dx
\\
=\,&\frac{1}{d_1d_2}\int_{\R^4}uU dx
-\frac{1}{d_1d_2}\int_{\R^4}U\pt_t wdx
-\frac{\lm_2}{d_1d_2}\int_{\R^4}Uwdx.
\eqntag
\label{eqn;uw}
\end{align*}
Recalling the Hardy--Littlewood--Sobolev inequality~\eqref{eqn;HLS-inequality},
\begin{align*}
\left|\iint_{\R^4\times\R^4}f(x)|x-y|^{-2}f(y)dxdy\right|
\le C_{HLS} \|f\|_{\frac43}^2,
\end{align*}
with optimal constant $C_{HLS}=(3/2)^{\frac12}\pi$,
we realize that
\begin{align*}
\frac{1}{d_1d_2}\int_{\R^4}uU dx
=\,&
\frac{1}{4\pi^2d_1d_2}\iint_{\R^4\times\R^4}u(x)|x-y|^{-2}u(y)dxdy
\\
\le\,&\frac{1}{4\pi^2d_1d_2}C_{HLS}\|u\|_{\frac43}^2.
\eqntag
\label{eqn;uw-estimate-HLS}
\end{align*}
In addition, it follows from~\eqref{eqn;FS-model} that $\pt_t W -d_2\Del W+\lm_2 W=U$ in $(0,\infty)\times\R^4$ since $W=E_4*w$. Thus,
\begin{align*}
-\frac{1}{d_1d_2}\int_{\R^4}U\pt_t wdx
=\,&
\frac{1}{d_1d_2}\int_{\R^4}(\pt_t W -d_2\Del W+\lm_2 W)\pt_t \Del W dx
\\
=\,&
-\frac{1}{d_1d_2}\|\N\pt_t W\|_2^2-\frac{1}{2d_1}\frac{d}{dt}\|\Del W\|_2^2
-\frac{\lm_2}{2d_1d_2}\frac{d}{dt}\|\N W\|_2^2
\\
=\,&
-\frac{1}{d_1d_2}\|\N\pt_t W\|_2^2-\frac{1}{2d_1}\frac{d}{dt}\|w\|_2^2
-\frac{\lm_2}{2d_1d_2}\frac{d}{dt}\|\N W\|_2^2.
\eqntag
\label{eqn;uw-estimate2}
\end{align*}
Similarly,
\begin{align*}
-\frac{\lm_2}{d_1d_2}\int_{\R^4}Uwdx
=\,&
\frac{\lm_2}{d_1d_2}\int_{\R^4}(\pt_t W -d_2\Del W+\lm_2 W)\Del Wdx
\\
=\,&-\frac{\lm_2}{2d_1d_2}\frac{d}{dt}\|\N W\|_2^2
-\frac{\lm_2}{d_1}\|\Del W\|_2^2
-\frac{\lm_2^2}{d_1d_2}\|\N W\|_2^2.
\eqntag
\label{eqn;uw-estimate3}
\end{align*}
Therefore we show from~\eqref{eqn;uw}--\eqref{eqn;uw-estimate3} that
\begin{align*}
\frac{1}{d_1}\int_{\R^4}uw dx
\le\,&
\frac{1}{4\pi^2d_1d_2}C_{HLS}\|u\|_{\frac43}^2
-\frac{1}{2d_1}\frac{d}{dt}\|w\|_2^2
-\frac{\lm_2}{d_1d_2}\frac{d}{dt}\|\N W\|_2^2
\\
&-\frac{1}{d_1d_2}\|\N\pt_t W\|_2^2
-\frac{\lm_2}{d_1}\|w\|_2^2
-\frac{\lm_2^2}{d_1d_2}\|\N W\|_2^2,
\end{align*}
and this concludes the proof of Lemma~\ref{lem;estimate-uw}.
\end{pr}

\vspace{5mm}

It is obvious from Lemma~\ref{lem;estimate-uw} that controlling the $L^1$-norm of $uw$ requires an estimate on the $L^{4/3}$-norm of $u$, which is provided by the next lemma, inspired by \cite[Lemma~2.3]{Na-Ya}.

\vspace{5mm}

\begin{lem}\label{lem;modified-sobolev}
	For any $\ep>0$ and $f\in W^{1,2}(\R^4)\cap L_+^1(\R^4)$, the following holds:
	\begin{align*}
	\|f\|_{\frac43}^2 \le C_S^2 (1+\ep)^2 \|f\|_1 \int_{\R^4}\frac{|\N f|^2}{1+f}dx+\frac{\sqrt{1+\ep}}{\ep}\left(\frac43\right)^{\frac32}\|f\|_1^{\frac32},
	\end{align*}
	where $C_S=2^{1/4}/(4\sqrt{\pi})$ is the optimal constant in the Sobolev inequality~\eqref{eqn;Sobolev-inequality-4/3}.
\end{lem}
\vspace{5mm}
\begin{pr}{Lemma~\ref{lem;modified-sobolev}}
For $f \in W^{1,2}(\R^4)\cap L_+^1(\R^4)$ and $\al>0$,
\begin{align}
	\int_{\R^4} f^{\frac43}dx=\,&\int_{ \{f>\al\}} (f-\al+\al)^{\frac43} dx +\int_{ \{f\le \al\}} f^{\frac43} dx \nonumber
	\\
	\le\,&\int_{ \{f>\al\}} \left((f-\al)^{\frac43}+\frac43\al f^{\frac13}\right)dx + \al^{\frac13}\int_{ \{f\le \al\}} f dx \nonumber
	\\
	\le\,& \int_{\R^4} (f-\alpha)_+^{\frac43}dx + \frac{4}{3} \alpha^{\frac13} \int_{\{ f>\alpha\}} f dx + \al^{\frac13}\int_{ \{f\le \al\}} f dx, \label{nd11}
\end{align}
where $(f-\al)_+:=\max\{f-\al,0\}$ and we have used the convexity inequality
\begin{equation*}
	(a+b)^{\frac43} \le a^{\frac43} +\frac43 b (a+b)^{\frac13},
	\quad
	(a,b)\in [0,\infty)^2.
\end{equation*}
Since
\begin{align*}
	\|\N (f-\al)_+\|_1=\,&
	\int_{ \{f>\al\}} |\N f| dx \le\, \left(\int_{ \{f>\al\}}(1+f)dx\right)^{\frac12}\left(\int_{ \{f>\al\}} \frac{|\N f|^2}{1+f} dx\right)^{\frac12}
	\\
	\le\,& \left(\frac{1+\al}{\al}\|f\|_1\right)^{\frac12} \left(\int_{ \R^4} \frac{|\N f|^2}{1+f} dx\right)^{\frac12},
\end{align*}
the function $(f-\alpha)_+$ belongs to $W^{1,1}(\R^4)$ and it follows from~\eqref{nd11} and the Sobolev inequality 
\begin{equation*}
\|f\|_{\frac43} \le C_S \|\N f\|_1,
\quad
f\in W^{1,1}(\R^4),
\eqntag
\label{eqn;best-sobolev}
\end{equation*}
with optimal constant $C_S=2^{1/4}/(4\sqrt{\pi})$ that
\begin{align*}
	\int_{\R^4} f^{\frac43}dx
	\le\,&
	\| (f-\al)_+\|_{\frac43}^{\frac43} + \frac43 \al^{\frac13} \|f\|_1
	\\
	\le\,&
	C_S^{\frac43} \|\N (f-\al)_+\|_1^{\frac43} + \frac43 \al^{\frac13} \|f\|_1.
\end{align*}
Collecting the above inequalities, we have for any $\ep>0$
	\begin{align*}
	\|f\|_{\frac43}^2
	\le\,&
	\left[
	C_S^{\frac43} \left(\frac{1+\al}{\al}\|f\|_1\right) ^{\frac23} \left(\int_{ \R^4} \frac{|\N f|^2}{1+f} dx\right)^{\frac23}
	+\frac43 \al^{\frac13} \|f\|_1
	\right]^{\frac32}
	\\
	\le\,&
	(1+\ep)C_S^2 \left(\frac{1+\al}{\al}\|f\|_1\right)\int_{ \R^4}\frac{|\N f|^2}{1+f}dx + \sqrt{\frac{1+\ep}{\ep}}
	\left(\frac43\right)^{\frac32} 
	\al^{\frac12} \|f\|_1^{\frac32}
	\end{align*}
	by applying 
	\begin{align*}
	(a+b)^{\frac32} \le (1+\ep)a^{\frac32}+ \sqrt{\frac{1+\ep}{\ep}}
	b^{\frac32},
	\quad
	(a,b)\in [0,\infty)^2.
	\end{align*}
	Then taking $\al=1/\ep$ leads us to
	\begin{align*}
	\|f\|_{\frac43}^2 
	\le
	(1+\ep)^2 C_S^2 \|f\|_1\int_{ \R^4}\frac{|\N f|^2}{1+f} dx+\frac{\sqrt{1+\ep}}{\ep}\left(\frac43\right)^{\frac32}\|f\|_1^{\frac32}.
	\end{align*}
	This concludes the proof.
\end{pr}

\vspace{5mm}
After this preparation, we are in a position to show the claimed boundedness of $\L(u,v,w)$ for solutions to~\eqref{eqn;FS-model} when $\|u_0\|_1$ satisfies the constraint~\eqref{eqn;bdd-assumption}.
\vspace{5mm}

\begin{pr}{Proposition~\ref{prop;bound-L}}
Recalling the definition~\eqref{eqn;L_0} and~\eqref{eqn;dissipative} of $\L$ and $\D_1$, 
Lemmas~\ref{lem;L_0} and~\ref{lem;estimate-uw} imply that
\begin{align*}
\frac{d}{dt}\L(u,v,w)+\D_1(u,v,w)
\le\,&
\frac{1}{4\pi^2d_1d_2}C_{HLS}\|u\|_{\frac43}^2
\\
&+\frac{1}{d_1}\|u_0\|_1^{\frac12}\|\pt_t v\|_2
+\frac{\lm_1}{d_1}\|u_0\|_1^{\frac12}\|v\|_2,
\end{align*}
from which we further obtain, after using~\eqref{eqn;mass} and Lemma~\ref{lem;modified-sobolev},
\begin{align*}
\frac{d}{dt}\L(u,v,w)+\D_1(u,v,w)
\le\,&
\frac{1}{4\pi^2d_1d_2}C_{HLS}C_S^2
(1+\ep)^2 \|u_0\|_1 \int_{\R^4}\frac{|\N u|^2}{1+u}dx
\\
&+\frac{\sqrt{1+\ep}}{\ep}C\|u_0\|_1^{\frac32}
+\frac{1}{d_1}\|u_0\|_1^{\frac12}\|\pt_t v\|_2
+\frac{\lm_1}{d_1}\|u_0\|_1^{\frac12}\|v\|_2
\end{align*}
for any $\ep>0$. Since
\begin{equation*}
\|u_0\|_1<\frac1{\sqrt{3}}(8\pi)^2d_1d_2
\end{equation*}
by~\eqref{eqn;bdd-assumption}, we may choose $\ep>0$ such that
\begin{align*}
\|u_0\|_1(1+\ep)^2=\frac12\left(\|u_0\|_1+\frac1{\sqrt{3}}(8\pi)^2d_1d_2\right)
\in \left(\|u_0\|_1,\frac{(8\pi)^2}{\sqrt{3}}d_1d_2\right).
\end{align*}
Then, noting that
\begin{align*}
\frac{C_{HLS} C_S^2}{4\pi^2}=\frac{\sqrt{3}}{(8\pi)^2}
\end{align*}
and
\begin{align*}
1-\delta := \frac{\sqrt{3}}{(8\pi)^2d_1d_2}
\frac12\left(\|u_0\|_1+\frac{(8\pi)^2d_1d_2}{\sqrt{3}}\right)=\frac12\left(1+\frac{\sqrt{3}}{(8\pi)^2d_1d_2}\|u_0\|_1\right) \in \left( \frac12 , 1 \right)
\end{align*}
we infer from Young's inequality, Lemma~\ref{lem;vw-bound}, and the positivity of $\lambda_1, \lambda_2$ that 
\begin{align*}
\frac{d}{dt}\L(u,v,w)+\D_1(u,v,w)
\le\,&
\frac{\sqrt{3}}{(8\pi)^2d_1d_2}
\frac12\left(\|u_0\|_1+\frac{(8\pi)^2d_1d_2}{\sqrt{3}}\right)
\int_{\R^4}\frac{|\N u|^2}{1+u}dx
\\
&+\frac{1}{d_1}\|u_0\|_1^{\frac12}\|\pt_t v\|_2
+\frac{\lm_1}{d_1}\|u_0\|_1^{\frac12}\|v\|_2+C
\\
\le\,&
(1-\del)\int_{\R^4}\frac{|\N u|^2}{1+u}dx
+ (1-\delta) \frac{\lm_1+\lm_2}{d_1} \|\pt_t v\|_2^2
+C \|v\|_2^2+C
\\
\le\,&(1-\del)\D_1(u,v,w) + C,
\end{align*}
which means that
 \begin{align*}
 \frac{d}{dt}\L(u,v,w)+\del\D_1(u,v,w)\le C.
 \eqntag
 \label{eqn;L-estimate}
 \end{align*}
 Next, by Lemma~\ref{lem;vw-bound}, the positivity of $\lm_1, \lm_2$, and~\eqref{eqn;FS-model},
\begin{align*}
\frac{1}{2d_1}\Ep_0(v)
=\,&\frac{1}{2d_1}\Big(\|\pt_t v\|_2^2
+d_1d_2\|\Del v\|_2^2
+(d_1\lm_2+d_2\lm_1)\|\N v\|_2^2
+\lm_1\lm_2\|v\|_2^2\Big)
\\
\le\,&
\frac{1}{2d_1}\|\pt_t v\|_2^2
+ \frac{d_2}{2d_1^1} \|\pt_tv+\lm_1 v-w\|_2^2+C
\\
\le\,&C\|\pt_t v\|_2^2+C\|w\|_2^2+C \le C (1+\D_1(u,v,w)).
\end{align*}
Also, since 
\begin{equation*}
	\|\nabla W\|_2^2 = \int_{\R^4} w W dx \le \frac{C_{HLS}}{4\pi^2} \|w\|_{\frac43}^2
\end{equation*}
by the Hardy--Littlewood--Sobolev inequality~\eqref{eqn;HLS-inequality} and $(1+r)\log{(1+r)} \le C \left( r + r^{\frac43} \right)$ for $r>0$, Lemma~\ref{lem;vw-bound},  Lemma~\ref{lem;modified-sobolev}, and~\eqref{eqn;mass} imply that
 \begin{align*}
\int_{ \R^4}(1+u)\log(1+u)dx
 +\frac{1}{2d_1}\|w\|_2^2 +&\, \frac{\lm_2}{d_1d_2}\| \N W \|_2^2
 \\
 \le\,&C\|u\|_{\frac43}^2+C(1+\D_1(u,v,w)) + C \|w\|_{\frac43}^2
 \\
 \le\,&C\int_{\R^4}\frac{|\N u|^2}{1+u}dx+C(1+\D_1(u,v,w))
 \\
 \le\,&C(1+\D_1(u,v,w)).
 \end{align*}
Consequently, gathering the previous inequalities, we find
 \begin{align*}
\L(u,v,w)\le C(1+\D_1(u,v,w) ).
 \end{align*}
 Together with~\eqref{eqn;L-estimate}, we end up with
 \begin{align*}
 \frac{d}{dt}\L(u,v,w)+\frac{\del}{2C}\L(u,v,w)+\frac{\del}{2}\D_1(u,v,w)\le \frac{\del}{2}+C,
 \end{align*}
from which Proposition~\ref{prop;bound-L} follows.
\end{pr}

\vspace{5mm}
\subsection{Further time-independent estimates}
We emphasize that the subsequent bounds differ from those used in the proofs of Proposition~\ref{prop;apriori-estimate} and Lemma~\ref{lem;u-1/2-log}~to establish global existence. 
\vspace{5mm}
\begin{lem}\label{lem;estimate-bdd-solution}
There is a constant $C>0$ such that
\begin{align*}
\int_t^{t+1} \left(\|u(s)\|_{\frac32}^{\frac32}+\|\N\pt_t v(s)\|_2^2+\|\N w(s)\|_2^2\right)ds\le C.
\end{align*}
\end{lem}
\vspace{5mm}

\begin{pr}{Lemma~\ref{lem;estimate-bdd-solution}}
It first follows from Proposition~\ref{prop;bound-L} that, for $t>0$,
\begin{align*}
\|(1+u(t))\log(1+u(t))\|_1+\|w(t)\|_2\le C,
\eqntag
\label{eqn;ulou-bdd}
\end{align*}
and
\begin{align*}
\int_t^{t+1}\left( \int_{ \R^4}\frac{|\N u(s)|^2}{1+u(s)}dx+ \|\N \pt_t v(s)\|_2^2 \right)ds\le C.
\eqntag
\label{eqn;dispative-bdd}
\end{align*}
We next infer from~\eqref{eqn;mass}, \eqref{eqn;ulou-bdd}, \eqref{eqn;dispative-bdd}, and Lemma~\ref{lem;L3/2-estimate} with $N=1$ that, for $t>0$,
\begin{align*}
\int_t^{t+1}\|u(s)\|_{\frac32}^{\frac32}ds \le\,&
 C \left(\sup_{t<s<t+1}\|(1+u(s)\log(1+u(s)))\|_1\right) \int_t^{t+1}\int_{ \R^4}\frac{|\N u(s)|^2}{1+u(s)}dxds+C\|u_0\|_1
\\
\le\,&C.
\eqntag
\label{eqn;L32-u-bdd}
\end{align*}
Finally, by~\eqref{eqn;FS-model}
\begin{align*}
\frac{d}{dt}\|w\|_2^2+2d_2\|\N w\|_2^2+2\lm_2\|w\|_2^2=2\int_{ \R^4}uw dx.
\end{align*}
Arguing as in the proof of Lemma~\ref{lem;u-1/2-log},
we use the H\"{o}lder, the Sobolev, and the Young inequalities, along with~\eqref{eqn;mass}, to obtain
\begin{align*}
2\int_{ \R^4}uwdx\le\,&2\|u\|_{\frac43}\|w\|_4 \le C \|u\|_1^{\frac14} \|u\|_{\frac32}^{\frac34} \|\N w\|_2 \le C \|u\|_{\frac32}^{\frac34} \|\N w\|_2
\\
\le\,&d_2\|\N w\|_2^2+ C\|u\|_{\frac32}^{\frac32}.
\end{align*}
Combining the above two inequalities gives
\begin{align*}
\frac{d}{dt}\|w\|_2^2+d_2\|\N w\|_2^2\le C \|u\|_{\frac32}^{\frac32}.
\end{align*}
Hence, after integrating on $(t,t+1)$ and using~\eqref{eqn;ulou-bdd} and~\eqref{eqn;L32-u-bdd},
\begin{align*}
d_2\int_t^{t+1}\|\N w(s)\|_2^2ds \le\|w(t)\|_2^2+C\int_t^{t+1}\|u(s)\|_{\frac32}^{\frac32}ds\le C,
\end{align*}
which completes the proof.
\end{pr}

\vspace{5mm}

Finally we proceed to show the proof of Theorem~\ref{thm;Bdd}. The argument used here is basically similar to that used in the proof of Theorem~\ref{thm;GE} to establish global existence. However,  it is necessary to carefully track the time dependence of the estimates and thus to modify the argument for the $L^p$-bound of $u$, because we only have the time-integrated uniform bounds on $\N\pt_t v$ and $\N w$ in $L^2\big((t,t+1)\times\R^4\big)$ stated in Lemma~\ref{lem;estimate-bdd-solution}.
\vspace{5mm}

\begin{pr}{Theorem~\ref{thm;Bdd}}
Arguing as in the proof of Theorem~\ref{thm;GE} with $p=3/2$, we obtain
\begin{equation*}
\frac{d}{dt}\|u(t)\|_{\frac32}^{\frac32} \le C\|u\|_{\frac32}^{\frac32}\left(\|\N w\|_2^2+\|\N \pt_t  v\|_{2}^2\right)
\end{equation*}
for $t>0$. Owing to Lemma~\ref{lem;estimate-bdd-solution}, we are in a position of applying the uniform Gronwall lemma (cf. Temam \cite[Lemma~III.1.1]{Te}) to conclude that
\begin{align*}
u\in L^{\infty}\big(0,\infty; L^{\frac32}(\R^4)\big).
\end{align*} 
We then infer from the parabolic regularity theory and the positivity of $\lm_1, \lm_2$ that
\begin{equation*}
w\in L^{\infty}\big(0,\infty; L^p(\R^4)\big),\quad p\in[1,6),
\end{equation*}
and 
\begin{equation*}
\N v \in L^{\infty}\big((0,\infty)\times \R^4\big).
\end{equation*}
Therefore we use the Nash--Moser iteration argument to obtain the uniform bound on the solution as in the proof of Theorem~\ref{thm;GE} and this concludes the proof of Theorem~\ref{thm;Bdd}.
\end{pr}
\vskip5mm
\noindent
{\bf Data availability statement.} 
All data that support the findings of this study are included within the article.

\vskip5mm
\noindent
{\bf Acknowledgments.} 
The work of the first author is partially supported by Grant-in-Aid for JSPS Fellows Grant Number 23KJ0150. The first author would like to thank Professor Takayoshi Ogawa for his great support. The authors would like to thank Professor Kentaro Fujie for his helpful suggestions and comments.

\end{document}